\documentclass[12pt]{amsart}

\usepackage{times,epsf,amssymb,amsmath,hyperref, rlepsf, color}

\begin{document}

\newtheorem{thm}{Theorem}[section]
\newtheorem{lem}[thm]{Lemma}
\newtheorem{cor}[thm]{Corollary}

\theoremstyle{definition}
\newtheorem{defn}[thm]{\bf{Definition}}

\theoremstyle{remark}
\newtheorem{rmk}[thm]{Remark}

\def\square{\hfill${\vcenter{\vbox{\hrule height.4pt \hbox{\vrule width.4pt height7pt \kern7pt \vrule width.4pt} \hrule height.4pt}}}$}

\newenvironment{pf}{{\it Proof:}\quad}{\square \vskip 12pt}

\title[$H$-Surfaces with Arbitrary Topology in $\BH^3$]{$H$-Surfaces with Arbitrary Topology in Hyperbolic 3-Space}
\author{Baris Coskunuzer}
\address{Koc University, Department of Mathematics, Istanbul 34450 Turkey}
\address{Max-Planck Institute for Mathematics, Bonn 53111 Germany}
\email{bcoskunuzer@ku.edu.tr}

\maketitle

%% User definitions:

\newcommand{\Si}{S^2_{\infty}({\Bbb H}^3)}
\newcommand{\SI}{S^n_{\infty}({\Bbb H}^{n+1})}
\newcommand{\PI}{\partial_{\infty}}

\newcommand{\BH}{\Bbb H}
\newcommand{\BHH}{{\Bbb H}^3}
\newcommand{\BR}{\Bbb R}
\newcommand{\BC}{\Bbb C}
\newcommand{\BZ}{\Bbb Z}

\newcommand{\e}{\epsilon}

\newcommand{\wh}{\widehat}

\newcommand{\A}{\mathcal{A}}
\newcommand{\C}{\mathcal{C}}
\newcommand{\p}{\mathcal{P}}
\newcommand{\R}{\mathcal{R}}
\newcommand{\B}{\mathcal{B}}
\newcommand{\h}{\mathcal{H}}
\newcommand{\T}{\mathfrak{T}}
\newcommand{\s}{\mathcal{S}}

\begin{abstract}

In this paper, we show that any open orientable surface $S$ can be properly embedded in $\BH^3$ as a minimizing $H$-surface for any $0\leq H<1$. We obtained this result by proving a version of the bridge principle at infinity for $H$-surfaces. We also show that any open orientable surface $S$ can be {\em nonproperly} embedded in $\BHH$ as a minimal surface, too.
\end{abstract}

\section{Introduction}

In this paper, we are interested in the existence of complete Constant Mean Curvature (CMC) surfaces in $\BHH$ of arbitrary topological type. CMC surfaces in the hyperbolic $3$-space has been an attractive topic for the last two decades. Especially after the substantial results on asymptotic Plateau problem, i.e. the existence and regularity of minimal surfaces in $\BHH$ by Anderson \cite{A1},\cite{A2}, and Hardt and Lin \cite{HL}, the generalizations of these results to CMC surfaces became interesting. In the following years, Tonegawa generalized Anderson's existence and Hardt and Lin's regularity results for CMC hypersurfaces \cite{To}.

Later, Oliveira and Soret studied the question of {\em "What kind of surfaces can be minimally embedded in $\BHH$?"} and showed that any finite topological type surface can be minimally embedded in $\BHH$ where the embedding is complete \cite{OS}. Then, Ros conjectured that any open surface (not necessarily finite topology) can be properly and minimally embedded in $\BHH$. Very recently, Francisco Martin and Brian White gave a positive answer to this conjecture, and showed that {\em any open orientable surface can be properly embedded in $\BHH$ as an area minimizing surface}  \cite{MW}. While in area minimizing case, there have been many great results on the realization of a surface of given topology in $\BHH$, there has been no result for CMC case so far in the literature.

In this paper, we address this problem, and generalize Martin and White's result to CMC surfaces ($H$-surfaces) for $0\leq H<1$. Our main result is as follows:

\begin{thm} Any open orientable surface can be properly embedded in $\BHH$ as a minimizing $H$-surface for \ $0\leq H<1$.
\end{thm}

In particular, this shows that any open orientable surface can be realized as a complete CMC surface with mean curvature $H$ in $\BHH$ where $0\leq H<1$. Also, $H=0$ case corresponds to the area minimizing case mentioned above \cite{MW}. While generalizing Martin and White's result to $H$-surfaces, we followed a similar but different path (See Final Remarks). In particular, the outline of the method is as follows.

Like \cite{MW}, we start with a simple exhaustion of the open orientable surface $S$ which is a decomposition into simpler surfaces $S_1\subset S_2 \subset ... S_n \subset ..$ where $S=\bigcup_{n=1}^\infty S_n$ \cite{FMM}. In other words, the surface $S$ can be constructed by starting with a disk $D=S_1$, and by adding $1$-handles iteratively, i.e. $S_{n+1}-int(S_n)$ is either a pair of pants attached to $S_n$ or a cylinder with a handle attached to $S_n$ (See Figure \ref{simpleexhaustion}). Hence after proving a version of bridge principle at infinity for $H$-surfaces, we started the construction with an $H$-plane in $\BHH$, say $S_1$. Then, if $S_{n+1}$ is a pair of pants attached to $S_n$, then we attach a bridge in $\Si$ to the corresponding component of $\PI S_n$. Similarly, if $S_{n+1}$ is a cylinder with a handle attached to $S_n$, then we attach two bridges successively to $\PI S_n$ (See Figures \ref{handle} and \ref{caps}). By iterating this process dictated by the simple exhaustion of $S$, we construct a properly embedded $H$ surface $\Sigma$ in $\BHH$ with the same topological type of $S$. 

After constructing properly embedded $H$-surfaces in $\BHH$, we turn to the question of {\em "What kind of surfaces can be nonproperly embedded in $\BHH$ as a minimal surface?"}. By {\em placing a bridge} between the nonproperly embedded minimal plane in $\BHH$ constructed in \cite{Co3}, and the minimal surface of desired topological type constructed above, we show that any open orientable surface can be minimally and nonproperly embedded in $\BHH$. 

\begin{thm} Any open orientable surface can be nonproperly embedded in $\BHH$ as a minimal surface.
\end{thm}

The organization of the paper is as follows. In the next section, we will give the basic definitions and results. In section 3, we will prove a version of bridge principle at infinity for $H$-surfaces in $\BHH$. In Section 4, we show the main result, the existence of properly embedded, complete minimizing $H$-surfaces in $\BHH$ of arbitrary topological type. In Section 5, we will show the existence of non-properly embedded minimal surfaces in $\BHH$ of arbitrary topological type. In section 6, we give some concluding remarks. Note that we postpone some technical steps to the appendix section at the end.

\subsection{Acknowledgements}

I would like to thank Brian White and Francisco Martin for very valuable conversations and remarks.

\section{Preliminaries}

In this section, we will overview the basic results which we use in the following sections. For further details, see \cite[Section 6]{Co1}

Let $\Sigma$ be a compact surface, bounding a domain $\Omega$ in some ambient Riemannian $3$-manifold. Let $A$ be the
area of $\Sigma$, and $V$ be the volume of $\Omega$. Let's vary $\Sigma$ through a one parameter family
$\Sigma_t$, with corresponding area $A(t)$ and volume $V(t)$. If $f$ is the normal component of the
variation, and $H$ is the mean curvature of $\Sigma$, then we get $A'(0) = -\int_\Sigma 2 H f$, and
$V'(0)=\int_\Sigma f$ where $H$ is the mean curvature.

Let $\Sigma$ be a surface with boundary $\alpha$. We fix a surface $M$ with $\partial M = \alpha$, and define $V(t)$ to be the volume of the domain bounded by $M$ and $\Sigma_t$. Now, we define a new functional as a combination of $A$ and $V$. Let $I_H(t)= A(t) + 2 H V(t)$. Note that $I_0(t)=A(t)$. If
$\Sigma$ is a critical point of the functional $I_H$ for any variation $f$, then this will imply $\Sigma$ has constant mean curvature $H$ [Gu]. Note that critical point of the functional $I_H$ is independent of the choice of the surface $M$ since if $\widehat{I}_H$ is the functional which is defined with a different surface $\widehat{M}$, then $I_H - \widehat{I}_H = C$ for some constant $C$. In particular, $H=0$ is the
special case of minimal surfaces and area minimizing surfaces, for which the theory is very well developed. We represent $H=0$ case in brackets [..] in the following definition. This definition describes well why CMC surfaces are considered as generalizations of minimal surfaces in a certain way.

%\begin{defn} $[H=0]$
%\begin{itemize}
%\item $\Sigma$ is a \textit{minimal surface} if it is critical point of $I_0$ (Area Functional) forany variation. Equivalently, $\Sigma$ has constant mean curvature $0$ at every point.
%\item A compact surface with boundary $\Sigma$ is an \textit{area-minimizing surface} if $\Sigma$ has the smallest area among surfaces with the same boundary (the absolute minimum of the functional $I_0$).
%\item A surface (not necessarily compact) is an \textit{area-minimizing surface} if any compact subsurface is an area minimizing surface.
%\item An area-minimizing surface $\Sigma$ with $\PI \Sigma = \Gamma$ is a {\em uniquely minimizing surface} if $\Gamma$ bounds a unique area minimizing surface in $\BHH$.
%\end{itemize}
%\end{defn}

%For general $H$, the generalization of minimal case to constant mean curvature $H$ case is as follows.

\begin{defn} {\bf i.} $\Sigma$ is called as \textit{$H$-surface} [minimal surface] if it is critical point of $I_H$ [$I_0$] for any variation.
Equivalently, $\Sigma$ has constant mean curvature $H$ ($0$) at every point.

{\bf ii.} A compact surface with boundary $\Sigma$ is a \textit{minimizing $H$-surface} [area minimizing surface] if $\Sigma$ is the absolute minimum of the functional $I_H$ [$I_0$] among surfaces with the same boundary.

{\bf iii.} A surface (not necessarily compact) is a \textit{minimizing $H$-surface} [area minimizing surface] if any compact subsurface is a minimizing $H$-surface [area minimizing surface].

{\bf iv.} A minimizing $H$-surface [area minimizing surface] $\Sigma$ with $\PI \Sigma = \Gamma$ is a {\em uniquely minimizing $H$-surface} [uniquely minimizing surface] if $\Gamma$ bounds a unique minimizing $H$-surface (area minimizing surface) in $\BHH$.
\end{defn}

\noindent \textbf{Notation:} From now on, we will call CMC surfaces with mean curvature $H$ as
\textit{$H$-surfaces} and we will assume $0\leq H<1$ unless otherwise stated. All the surfaces are assumed to be orientable unless otherwise stated.

Now, we will give the basic results on $H$-surfaces in hyperbolic space. The following existence result is given in any dimension.

\begin{lem} \label{existence} \cite{To}, \cite{AR}
Let $\Gamma$ be a codimension-$1$ closed submanifold in $\SI$, and let $|H|<1$. Then there exists a
minimizing $H$-hypersurface $\Sigma^n$ in $\BH^{n+1}$ where $\PI \Sigma^n = \Gamma$. Moreover,
any such $H$-hypersurface is smooth outside of a closed set of Hausdorff dimension $n-7$.
\end{lem}

Beside the existence results,  Tonegawa studied the regularity at infinity in \cite{To}, and obtained the following result.

\begin{thm} \label{bounreg} \cite{To}
Let $\Gamma$ be a collection of $C^\infty$-smooth disjoint simple closed curves in $\Si$. Let $\Sigma$ be an $H$-surface in $\BH^3$ with $\PI \Sigma = \Gamma$. Then, $\Sigma \cup \Gamma$ is a $C^\infty$ submanifold with boundary in $\overline{\BHH}$.
\end{thm}

Note also that by using some barrier arguments, it is not hard to show that if $\theta_H$ is the intersection angle at infinity between an $H$-surface and the asymptotic boundary $\Si$, then $\cos{\theta_H}=H$ \cite{To}.

%We should also note that Nelli and Spruck showed existence of a $H$-hypersurface asymptotic to $C^{2,\Gamma}$ codimension-1 submanifold $\alpha$ which is the boundary of a mean convex domain in $\SI$ by using analytic techniques in \cite{NS}. Later, Guan and Spruck generalized this result to $C^{1,1}$ codimension-1 submanifolds bounding star shaped domains in $\SI$ \cite{GS}.

%\begin{lem} \cite{GS} Let $\Omega$ be a star shaped (mean convex in \cite{NS}) domain in $\SI$ where $\alpha=\partial \Omega$ is $C^{1,1}$ ($C^{2,\Gamma}$ in \cite{NS}) codimension-1 submanifold in $\SI$. Then, for any $0<H<1$, there exists a complete smoothly embedded CMC surface $\Sigma$ with mean curvature $H$ and $\PI \Sigma = \alpha$. Moreover, $\Sigma$ can be represented as a graph of a function $u \in C^{1,1}(\overline{\Omega})$ ($u \in C^{2,\Gamma}(\overline{\Omega})$ in \cite{NS}).\end{lem}

%Even though this second existence result is for fairly restricted class of asymptotic boundary data (star shaped condition), the CMC surfaces obtained are smoothly embedded with no singularity in any dimension (unlike the first one), and they can be represented as a graph like $x_{n+1} = u$ for a function $u \in C^{1,1}(\overline{\Omega})$ in half space model for $\BHH$.

The following fact is known as maximum principle.

\begin{lem} \label{maxprinciple}
Let $\Sigma_1$ and $\Sigma_2$ be two surfaces in a Riemannian manifold which intersect at a common point
tangentially. If $\Sigma_2$ lies in positive side (mean curvature vector direction) of $\Sigma_1$ around the common
point, then $H_1$ is less than or equal to $H_2$ ($H_1 \leq H_2$) where $H_i$ is the mean curvature of $\Sigma_i$ at
the common point. If they do not coincide in a neighborhood of the common point, then $H_1$ is strictly less than $H_2$ ($H_1<H_2$).
\end{lem}

%For the number of solutions to the asymptotic Plateau problem, by using analytic techniques, Nelli and Spruck %generalized Anderson's uniqueness result for mean convex domains in area minimizing surfaces case to CMC context in %\cite{NS}. Then, Guan and Spruck extended Hardt and Lin's uniqueness results for star-shaped domains in area %minimizing surfaces case to CMC surfaces in hyperbolic space in \cite{GS}.

%\begin{lem} \label{uniqstar} \cite{GS}
%Let $\Omega$ be a star shaped (mean convex in \cite{NS}) domain in $\SI$ where $\alpha=\partial
%\Omega$ is $C^{1,1}$ ($C^{2,\Gamma}$ in \cite{NS}) codimension-1 submanifold in $\SI$. Then, for
%any $0\leq H<1$, there exists a unique complete CMC surface $\Sigma$ with mean curvature $H$
%and $\PI \Sigma = \alpha$.
%\end{lem}

Now, we will quote the following result from \cite{Co2}, which is used to prove the genericity of uniquely minimizing $H$-surfaces in $\BHH$. This lemma will also be an important tool for us to prove the bridge principle at infinity.

%\begin{lem} \label{genuniq} \cite[Theorem 4.1]{Co2}Let $A$ be the space of simple closed curves in $\Si$, and let $A'\subset A$ be the subspace containing the simple closed curves bounding a unique minimizing $H$-surface in $\BHH$. Then $A'$ is generic in $A$ in Baire sense. \end{lem} Next, we will quote another result from \cite{Co2}.

\begin{lem} \label{canonical}\cite[Lemma 4.1]{Co2}
Let $\Gamma$ be a collection of simple closed curves in $\Si$. Then either there exists a unique minimizing
$H$-surface $\Sigma$ in $\BHH$ with $\PI \Sigma = \Gamma$, or there are two canonical disjoint extremal
minimizing $H$-surfaces $\Sigma^+$ and $\Sigma^-$ in $\BHH$ with $\PI \Sigma^\pm = \Gamma$.
\end{lem}

Now, we will show that if two disjoint collection of simple closed curves, say $\beta_1$ and $\beta_2$, does not "link" each other in $\Si$, then the minimizing $H$-surfaces $T_1$ and $T_2$ in $\BHH$ with $\PI T_i=\beta_i$ must be disjoint.

\begin{lem} \label{disjoint} Let $\Omega_1$ and $\Omega_2$ be two open subsets (not necessarily connected) in $\Si$ with $\overline{\Omega}_1\cap\overline{\Omega}_2=\emptyset$. Let $\beta_i=\partial \overline{\Omega}_i$ be smooth curves. Then, if $T_1$ and $T_2$ are two minimizing $H$-surfaces in $\BHH$ with $\PI T_i=\beta_i$, then $T_1\cap T_2=\emptyset$.
\end{lem}

\begin{pf} We will basically adapt the technique in \cite[Theorem 3.2]{Co2} to this case. Since $\beta_i$ is smooth, $T_i\cup\beta_i$ is smoothly embedded in $\overline{\BHH}$ by Lemma \ref{bounreg}. Since $T_i$ is connected, $T_i$ separates $\BHH$ into two regions, say $\BHH-T_i=\Delta^+_i\cup\Delta^-_i$ where $\PI \Delta^+_i = \Omega_i$.

Assume that $T_1\cap T_2\neq \emptyset$. Then by maximum principle (Lemma \ref{maxprinciple}), $\Delta^+_1 \cap \Delta^+_2 \neq \emptyset$. Since $\overline{\Omega}_1\cap \overline{\Omega}_2=\emptyset$, then $W=\Delta^+_1 \cap\Delta^+_2$ is in the compact part of $\BHH$. Let $\partial \overline{W}\cap T_1 = \Sigma_1$ and $\partial \overline{W}\cap T_2 = \Sigma_2$. Then, $\partial \Sigma_1 = \partial \Sigma_2 = T_1 \cap T_2$ is a collection of simple closed curves, say $\tau$.

Now, recall that $T_1$ and $T_2$ are both minimizing $H$-surfaces in $\BHH$, and hence, the compact subsurfaces $\Sigma_1$ and $\Sigma_2$ are minimizing $H$-surfaces with the same boundary $\tau$, i.e. $I_H(\Sigma_1)=I_H(\Sigma_2)$. Now, let $S_1$ be a compact subsurface of $T_1$ with $\Sigma_1 \subset int(S_1)$. Define $S_1'=(S_1 -\Sigma_1)\cup \Sigma_2$. Clearly, $\partial S_1 = \partial S_1'$. Since $I_H(\Sigma_1)=I_H(\Sigma_2)$, then $I_H(S_1)=I_H(S_1')$ by construction. As $S_1$ is minimizing $H$-surface, and $I_H(S_1)=I_H(S_1')$, then $S_1'$ is also a minimizing $H$-surface with the same boundary. However, $S_1'$ has codimension-$1$ singularity along $\tau$. This contradicts to the regularity theorem for minimizing $H$-surfaces.
\end{pf}

\begin{rmk} \label{disj} Note that in the lemma above, we can take the open subsets $\Omega_1$ and $\Omega_2$ in $\Si$ with the condition $\overline{\Omega}_1\subset \Omega_2$. This is because this would be equivalent to $\overline{\Omega}_1\cap \overline{int(\Omega_2^c)}=\emptyset$ which satisfies the assumption in the lemma (i.e. take the second open set $\Omega_2'$ as $int(\Omega_2^c)$). In particular, this is a nonlinking condition for $\beta_1$ and $\beta_2$ in $\Si$.
\end{rmk}

We will finish this section with the following definition.

\begin{defn} Let $S$ be a noncompact surface. An embedding $\varphi:S\to X$ is {\em proper} if for any compact subset $K$ of $X$, $\varphi^{-1}(K)$ is compact in $S$. A surface $\Sigma$ is {\em properly embedded} in $X$ if there exists a proper embedding $\varphi:S\to X$ with $\varphi(S)=\Sigma$. Equivalently, $\Sigma$ is {\em properly embedded} in $X$ if $\overline{\Sigma}=\Sigma$ where $\overline{\Sigma}$ is the closure of $\Sigma$ in $X$.
\end{defn}

\section{A Bridge Principle at Infinity for H-surfaces}

\begin{figure}[h]
\begin{center}
$\begin{array}{c@{\hspace{.2in}}c}

\relabelbox  {\epsfxsize=2.5in \epsfbox{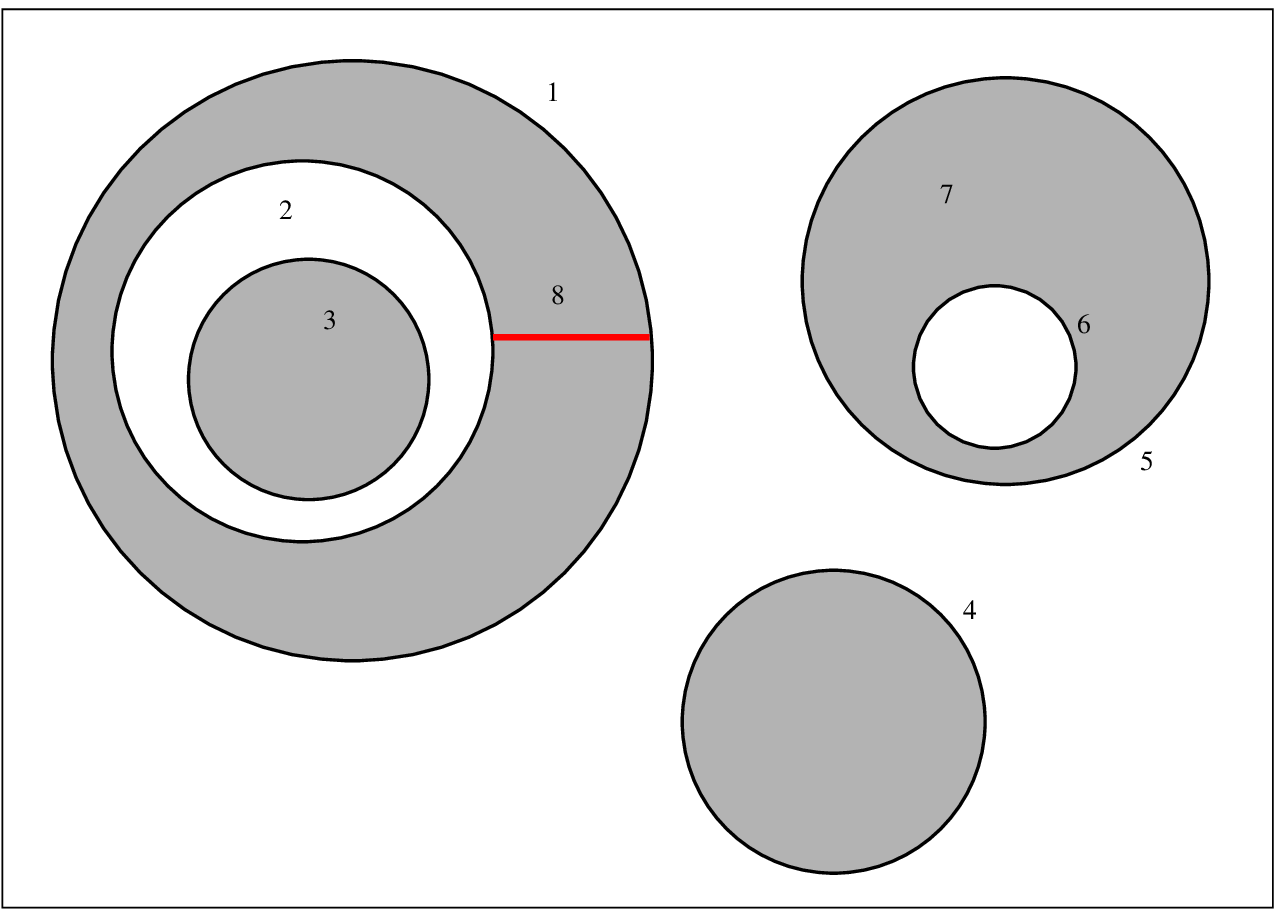}} \relabel{1}{$\gamma_1$} \relabel{2}{$\gamma_2$} \relabel{3}{$\gamma_3$} \relabel{4}{$\gamma_4$} \relabel{5}{$\gamma_5$} \relabel{6}{$\gamma_6$} \relabel{7}{$X^+$} \relabel{8}{$\alpha$}  \endrelabelbox &

\relabelbox  {\epsfxsize=2.5in \epsfbox{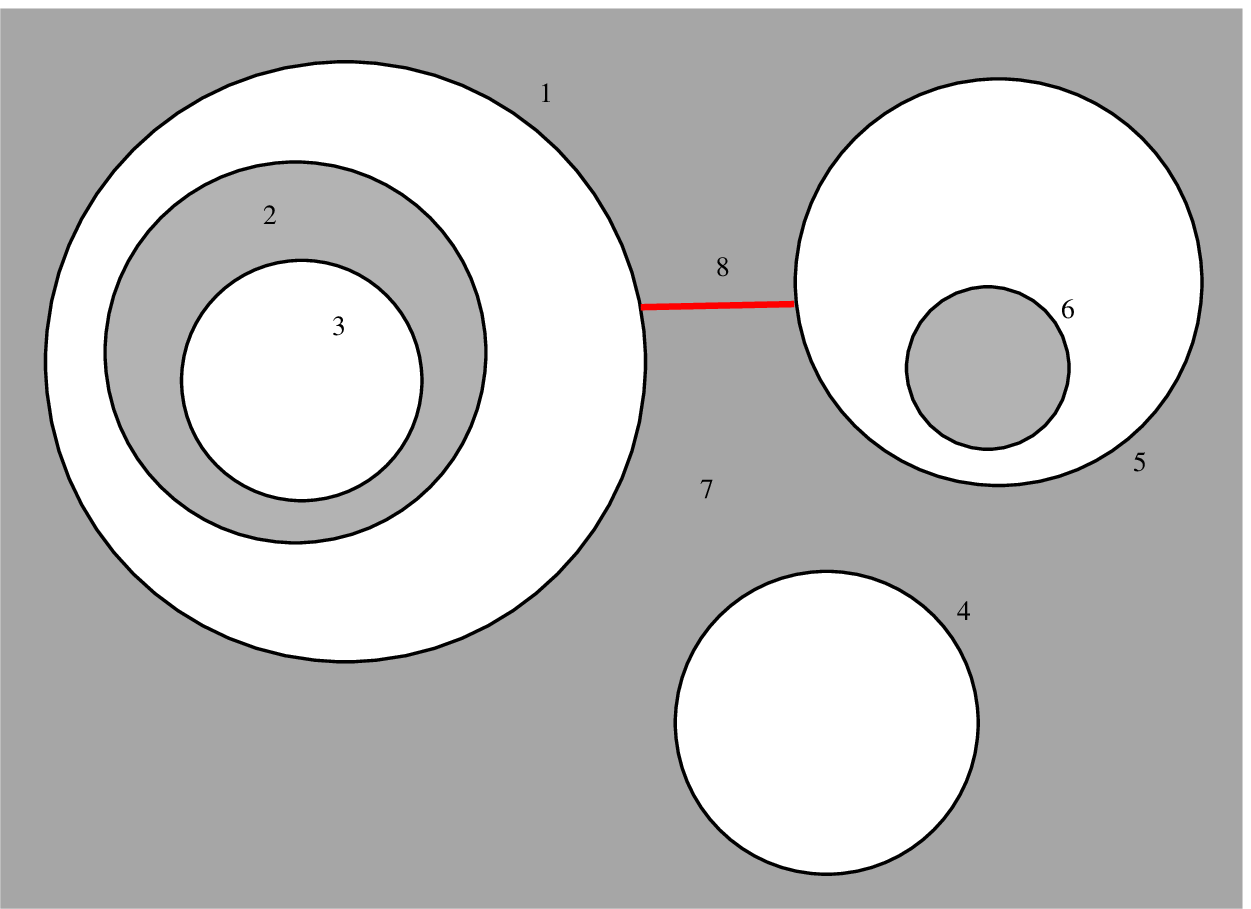}} \relabel{1}{$\gamma_1$} \relabel{2}{$\gamma_2$} \relabel{3}{$\gamma_3$} \relabel{4}{$\gamma_4$} \relabel{5}{$\gamma_5$} \relabel{6}{$\gamma_6$} \relabel{7}{$X^+$} \relabel{8}{$\alpha$} \endrelabelbox \\ [0.4cm]
\end{array}$

\end{center}
\caption{\label{XX} \small Let $\Gamma=\gamma_1\cup \gamma_2\cup...\cup\gamma_6$ be a collection of simple closed curves in $\Si$. Then, $\Gamma$ separates $\Si$ into two parts, say $X^+$ and $X^-$ with $\partial X^+=\partial X^-=\Gamma$. We call the gray region which contains the bridge $\alpha$ as $X^+$. In the pictures left and right, the situations are given when the bridge $\alpha$ is in different sides of $\Gamma$. }
\end{figure}

In this part, we will generalize the bridge principle at infinity by Martin and White to $H$-surfaces in $\BHH$ by using different techniques (See Final Remarks). Let $\Gamma$ be a finite collection of smooth simple closed curves in $\Si$. Let $\alpha$ be a smooth arc in $\Si$ which meets $\Gamma$ orthogonally, and satisfying $\Gamma \cap \alpha = \partial \alpha$. $\Gamma$ separates $\Si$ into two regions, say $\Si-\Gamma=X^+\cup X^-$ with $\partial X^+=\partial X^-=\Gamma$ (Notice that if $U\subset S^2$ and $\partial U= \gamma$ then $\partial U^c = \gamma$, too.). Of course, $X^+$ or $X^-$ may have more than one component as $\Gamma$ may not be connected. Let $X^+$ be the region which contains $\alpha$. See Figure \ref{XX}.

\begin{figure}[b]
\begin{center}
$\begin{array}{c@{\hspace{.2in}}c}

\relabelbox  {\epsfxsize=2.5in \epsfbox{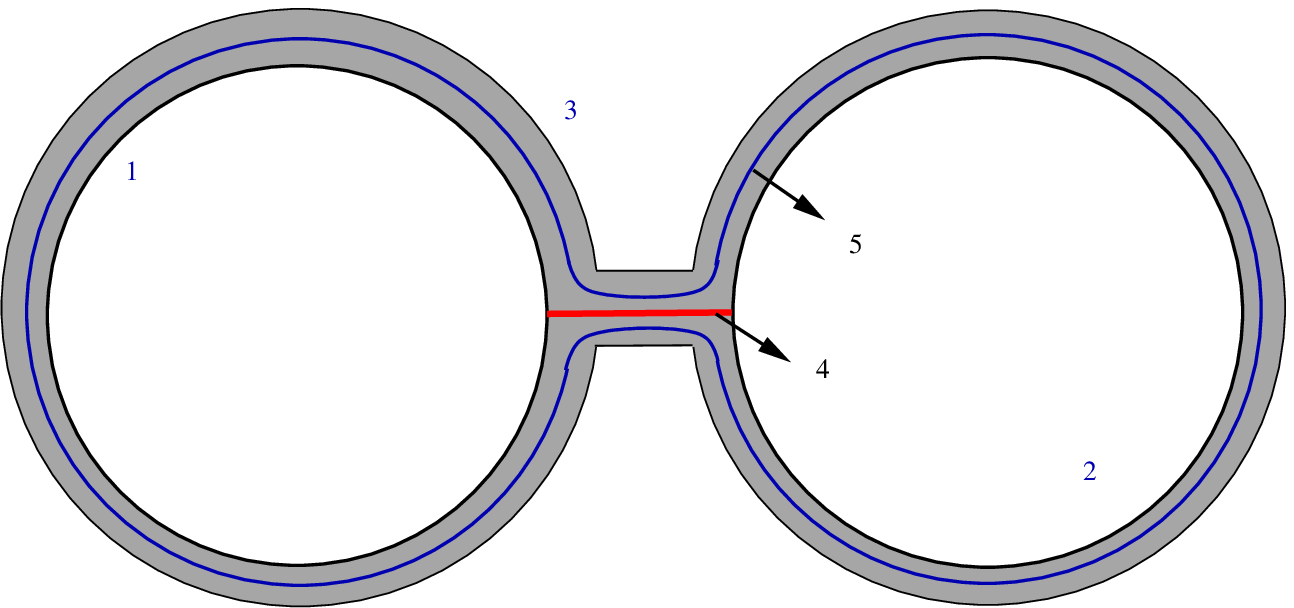}} \relabel{1}{$\gamma_1$} \relabel{2}{$\gamma_2$} \relabel{3}{$\Gamma_\delta$} \relabel{4}{$\alpha$} \relabel{5}{$\Gamma_t$} \endrelabelbox &

\relabelbox  {\epsfxsize=2.5in \epsfbox{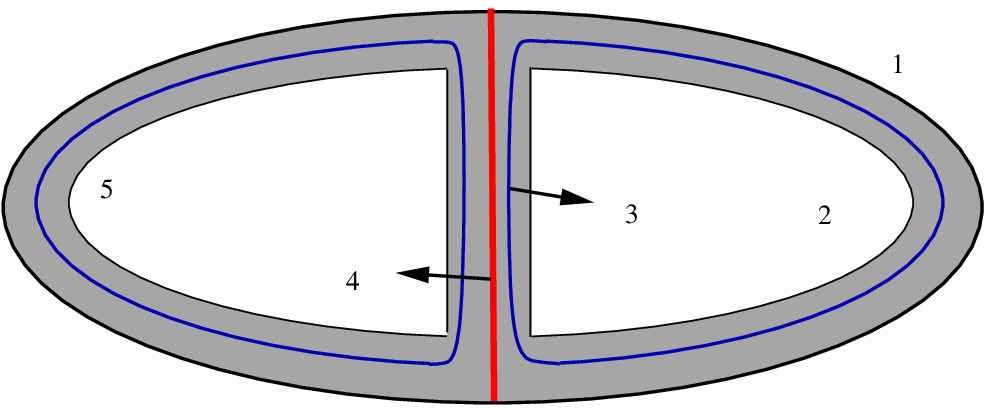}} \relabel{1}{$\Gamma$} \relabel{2}{$\Gamma^1_\delta$} \relabel{3}{$\Gamma_t$} \relabel{4}{$\alpha$} \relabel{5}{$\Gamma^2_\delta$}  \endrelabelbox \\ [0.4cm]
\end{array}$

\end{center}
\caption{\label{foliation} \small The gray region represents an open planar region $Y$ with $\partial Y=\Gamma\cup\Gamma_\delta$. $\{\Gamma_t \ | \ 0<t<\delta\}$ is a foliation of $Y$ by smooth curves. The red arc $\alpha$ is the bridge. In the left, $\alpha$ is connecting different components $\gamma_1$ and $\gamma_2$ of $\Gamma$, and  hence $\Gamma_t$ is connected for $0<t\leq \delta$. In the right, $\Gamma$ is connected and $\Gamma_t$ has 2 components for $0<t\leq \delta$. }
\end{figure}

Let $N_\epsilon(\Gamma)$ be the $\epsilon$ neighborhood of $\Gamma$ in $\Si$. Let $N^+(\Gamma)=X^+\cap N_\epsilon(\Gamma)$. Let $N_\epsilon(\alpha)$ be the $\epsilon$ neighborhood of $\alpha$ in $\Si$. Let $N^+_\e(\alpha)=X^+\cap N_\epsilon(\alpha)$. In other words, $N^+_\e(\alpha)$ is the component of $N_\e(\alpha)-\Gamma$ containing $\alpha$. Let $Y$ be the open planar region $N^+(\Gamma)\cup N^+(\alpha)$ in $X^+$ where $\partial \overline{Y} = \Gamma \cup \Gamma_\delta$. Foliate $Y$ by smooth collection of simple closed curves $\{\Gamma_t \ | \ 0<t<\delta\}$ such that $\Gamma_t\to\Gamma\cup\alpha$ as $t\searrow 0$, and $\Gamma_t\to \Gamma_\delta$ as $t\nearrow \delta$ (See Figure \ref{foliation}).

Notice that if the endpoints of $\alpha$ are in the same component of $\Gamma$, then $\sharp(\Gamma_t)=\sharp(\Gamma) +1$ and if the endpoints of $\alpha$ are in different components of $\Gamma$, then $\sharp(\Gamma_t)=\sharp(\Gamma) -1$ where $\sharp(\Gamma)$ represents the number of components of $\Gamma$.

Now, we prove a bridge principle at infinity for $H$-surfaces in $\BHH$. Note that we postpone some technical steps to the appendix section, which show that there is no genus on the bridge near infinity by using $H$-strips and $H$-skillets technique of \cite{MW}.

%In particular, we will show that if $S$ is uniquely minimizing $H$-surface with $\PI S = \Gamma$, then there exists sufficiently small $t>0$ such that the $H$-surface $S_t$ with $\PI S_t=\Gamma_t$ is uniquely minimizing $H$-surface with the desired topology, i.e. $S_t \simeq S\cup \widehat{N}_\epsilon(\alpha)$ where $\widehat{N}_\epsilon(\alpha)$ is the $\epsilon$-strip (bridge) along $\alpha$ in $\Si$, i.e. $\widehat{N}_\epsilon(\alpha)= \overline{X^+}\cap N_\epsilon(\alpha)$.

\begin{thm} \label{bridge} [A Bridge Principle at Infinity]

Let $S$ be a properly embedded, uniquely minimizing connected $H$-surface in $\BHH$ where $\PI S=\Gamma$ is a finite collection of disjoint smooth curves. Assume also that $\overline{S}$ has finite genus. Let $\alpha$ be a smooth arc in $\Si$ with $\Gamma\cap \alpha = \partial \alpha$ and $\Gamma\perp \alpha$. Consider the family of curves $\{\Gamma_t \ | \ 0<t<\delta\}$ constructed above. Then there exists a sufficiently small $t>0$ such that $\Gamma_t$ bounds a unique minimizing $H$-surface $S_t$ where $S_{t}$ is homeomorphic to $S\cup N^+_\epsilon(\alpha)$.
\end{thm}

\begin{pf} First, by Lemma \ref{existence}, for any $\Gamma_t\subset \Si$, there exists a minimizing $H$-surface $S_t$ with $\PI S_t = \Gamma_t$.

\vspace{.2cm}

\noindent {\bf Step 1:} For sufficiently small $t>0$, $S_t \simeq S\cup N^+_\epsilon(\alpha)$.

\vspace{.2cm}

\begin{pf} In this step, we will mainly use the techniques of \cite{MW}. As $t_n\searrow 0$, $S_{t_n}\to T$ where $T$ is a minimizing $H$-surface in $\BHH$ with $\PI T\subset \Gamma\cup \alpha$. By using the linking argument in \cite{MW}, one can show that $\PI T = \Gamma$. Since $S$ is uniquely minimizing $H$-surface with $\PI S=\Gamma$, $S=T$. Hence $S_{t_n}\to S$ and the convergence is smooth on compact sets.

We will use the upper half space model for $\BHH$. Assume that for $\e_n\searrow 0$, there exists $0<t_n<\e_n$ such that $S_{t_n}$, say $S_n$ for short, is not homeomorphic to $\wh{S}=S\cup N^+_{\epsilon}(\alpha)$. Since the number of boundary components are same, this means $S_n$ and $\wh{S}$ have different genus.

%Let $T$ be a minimizing $H$-surface in $\BHH$ with smooth asymptotic boundary $\gamma$, i.e. $\PI T = \gamma$.

Let $\R_a=\{0\leq z \leq a\}$ in $\overline{\BHH}$. In the appendix section, we show that there exists $a_\Gamma>0$ such that for sufficiently large $n$, $S_n\cap R_{a_\Gamma}$ has no genus.

Now, let $\mathcal{K}_a=\{z\geq a\}$ and let $S^a=S\cap \mathcal{K}_a$. Then, since $S_n\to S$ converge smoothly on compact sets,  $S^a_{t_n}\to S^a$ smoothly. Hence, by Gauss-Bonnet, $S^a_{t_n}$ and $S^a$ must have same genus. By above, this implies for sufficiently large $n$, $S_{t_n}$ and $S$ must have the same genus. This is a contradiction.

Hence, for sufficiently small $\epsilon'>0$, we will assume that for $0<t<\epsilon'$, $S_t$ is homeomorphic to $S\cup N^+_\epsilon(\alpha)$.
\end{pf}

\noindent {\bf Step 2:} For all but countably many $0<t<\epsilon'$, $\Gamma_t$ bounds a unique minimizing $H$-surface in $\BHH$.

\vspace{.2cm}

\begin{pf} In this step, we mainly use the techniques from \cite[Theorem 4.1]{Co2}. By Lemma \ref{disjoint} (see also Remark \ref{disj}), for any $0<t_1<t_2<\epsilon'$, if $S_1$ and $S_2$ are minimizing $H$-surfaces with $\PI S_i = \Gamma_{t_i}$, then $S_1$ and $S_2$ are disjoint. By Lemma \ref{canonical}, if $\Gamma_s$ does not bound unique minimizing $H$-surface, then we can define two disjoint canonical minimizing $H$-surfaces $S^+_s$ and $S^-_s$ with $\PI S^\pm_s = \Gamma_s$. Hence, $S_s^+\cup S_s^-$ separates a region $V_s$ from $\BHH$. If $\Gamma_s$ bounds a unique minimizing $H$-surface $S_s$, then let $V_s=S_s$. Notice that by lemma \ref{disjoint}, $S_t\cap S_s=\emptyset$ for $t\neq s$, and hence $V_t\cap V_s=\emptyset$ for $t\neq s$.

%In other words, by taking a sequence $t^+_n\searrow s$ and $t^-_n\nearrow s$, we get two limit minimizing $H$-surfaces as $S_{t^+_n} \to S_s^+$ and $S_{t^-_n} \to S_s^-$. By construction, if $S_s^+=S_s^-$ then $\Gamma_s$ bounds a unique minimizing $H$-surface. If not, then $S_s^+\cap S^-_s =\emptyset$ as $S_{t^+_n}\cap S_{t^-_m}=\emptyset$ by Lemma \ref{disjoint} for any $n, m$.  Moreover, this region is canonical, and independent of the choice of the sequences \cite{Co2}.

Now, consider a short arc segment $\eta$ in $\BHH$ with one endpoint is in $S_{t_1}$ and the other end point is in $S_{t_2}$ where $0<t_1<t_2<\epsilon'$. Hence, $\eta$ intersects all minimizing $H$-surfaces $S_t$ with $\PI S_t = \Gamma_t$ where $t_1\leq t \leq t_2$. Now for $t_1<s<t_2$, define the {\em thickness} $\lambda_s$ of $V_s$ as $\lambda_s=|\eta\cap V_s|$, i.e. $\lambda_s$ is the length of the piece of $\eta$ in $V_s$.Hence, if $\Gamma_s$ bounds more than one $H$-surface, then the thickness is not $0$. In other words, if $\lambda_s=0$, then $\Gamma_s$ bounds a unique $H$-surface in $\BHH$.

As $V_t\cap V_s=\emptyset$ for $t\neq s$, $\sum_{t_1}^{t_2}\lambda_s<|\eta|$. Hence, for only countably many $s\in[t_1,t_2]$, $\lambda_s>0$. This implies for all but countably many $s\in[t_1,t_2]$, $\lambda_s=0$, and hence $\Gamma_s$ bounds a unique minimizing $H$-surface. Similarly, this implies for all but countably many $s\in[0,\epsilon']$, $\Gamma_s$ bounds a unique $H$-surface. The proof follows.
\end{pf}

Steps 1 and 2 implies the existence of smooth curve $\Gamma_t$ with $0<t<\e'$ for any $\e'$, where $\Gamma_t$ bounds a unique minimizing $H$-surface $S_t$, and $S_t$ has the desired topology, i.e. $S_t\simeq S\cup \widehat{N}_{\epsilon}(\alpha)$.
\end{pf}

\begin{figure}[b]

\relabelbox  {\epsfxsize=3.5in

\centerline{\epsfbox{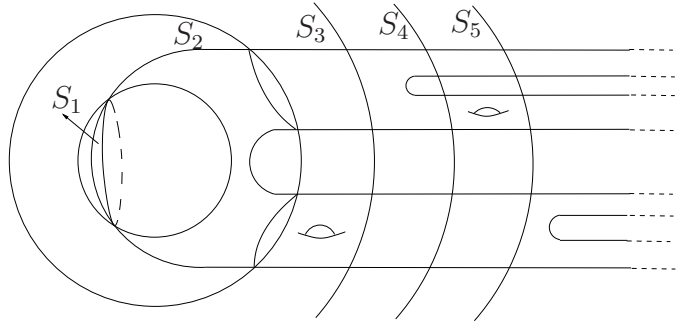}}}

\relabel{1}{$S_1$}

\relabel{2}{$S_2$}

\relabel{3}{$S_3$}

\relabel{4}{$S_4$}

\relabel{5}{$S_5$}

\endrelabelbox

\caption{\label{simpleexhaustion} \small In the simple exhaustion of $S$, $S_1$ is a disk, and $S_{n+1}-S_n$ contains a unique nonannular part, which is a pair of pants (e.g. $S_4-S_3$), or a cylinder with a handle (e.g. $S_3-S_2$). }

\end{figure}

\section{The Construction of Properly Embedded $H$-surfaces}

Now, we are going to prove the main existence result for properly embedded minimizing $H$-surfaces in $\BHH$ with arbitrary topology. In this part, we will mainly follow the techniques in \cite{MW}. In particular, for a given surface $S$, we will start with a compact exhaustion of $S$, $S_1\subset S_2\subset ... S_n \subset ...$, and by using the bridge principle in the previous section, we construct the minimizing $H$-surface with the desired topology.

In particular, by \cite{FMM}, for any open orientable surface $S$, there exists a simple exhaustion. A simple exhaustion $S_1\subset S_2\subset ... S_n \subset ...$ is the compact exhaustion with the following properties: $S_1$ is a disk, and $S_{n+1}-S_n$ would contain a {\em unique nonannular piece} which is either a cylinder with a handle (a torus with two holes), or a pair of pants by \cite{FMM} (See Figure \ref{simpleexhaustion}).

Hence, by starting with a round circle in $\Si$ which bounds a unique $H$-surface in $\BHH$ (a spherical cap), adding the bridges dictated by the simple exhaustion, we get a minimizing $H$-surface with the desired topology.

\begin{thm} \label{main} Any open orientable surface $S$ can be embedded in $\BHH$ as a minimizing $H$-surface $\Sigma$.
\end{thm}

\begin{pf} Let $S$ be an open orientable surface. Now, we inductively construct the minimizing $H$-surface $\Sigma$ in $\BHH$ which is diffeomorphic to $S$. Let $S_1\subset S_2\subset ... S_n \subset ...$ be a simple exhaustion of $S$, i.e. $S_{n+1}-S_n$ contains a unique nonannular piece which is either a cylinder with a handle, or a pair of pants.

\begin{figure}[b]
\begin{center}
$\begin{array}{c@{\hspace{.3in}}c}

\relabelbox  {\epsfysize=1.7in \epsfbox{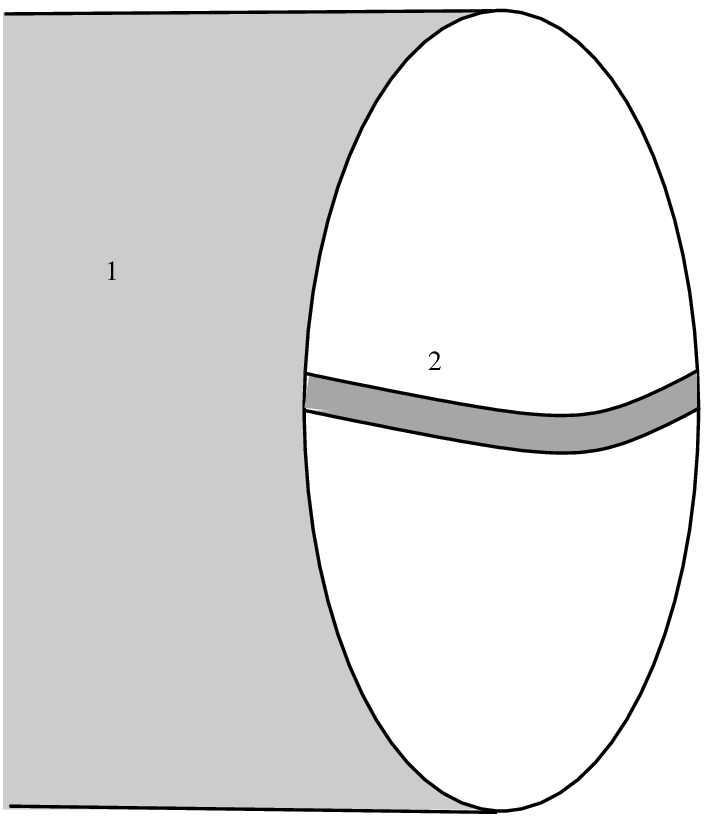}} \relabel{1}{$S_n$} \relabel{2}{$\B_n$}  \endrelabelbox &

\relabelbox  {\epsfysize=1.7in \epsfbox{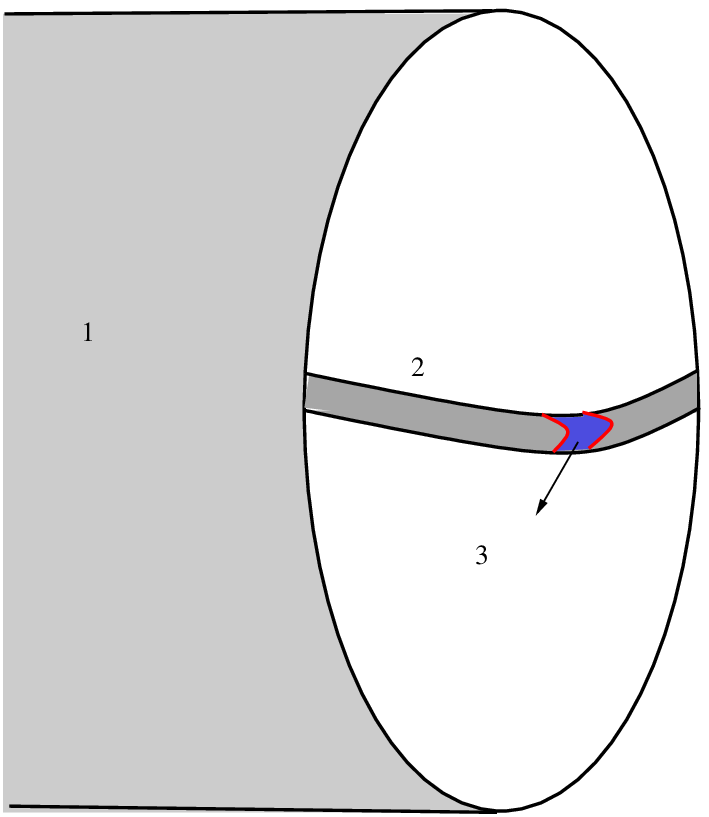}} \relabel{1}{$S_n$} \relabel{2}{$\B_n$} \relabel{3}{$\wh{\B}_n$}   \endrelabelbox \\
\end{array}$

\end{center}
\caption{\label{handle} \small If $S_{n+1}-S_n$ contains a pair of pants in the simple exhaustion, we add a bridge $\B_n$ so that $S_n\cup\B_n\simeq S_{n+1}$ (left). If  $S_{n+1}-S_n$ contains a cylinder with a handle, then we add a handle $\h_n$ so that $S_n\cup\h_n\simeq S_{n+1}$. Here the handle $\h_n$ is just successive two bridges, i.e $\h_n=\B_n\cup\wh{\B}_n$ (right). }
\end{figure}

Here, adding a bridge to the same boundary component of a surface would correspond to the attaching a pair of pants. Adding two bridges successively to the same boundary component would correspond to the attaching a cylinder with a handle. In particular, if $\C$ is the boundary component in $\partial S_n$ and the annulus $\A$ is a small neighborhood of $\C$ in $S_n$, then $\A\cup \B_n$ would be a pair of pants, where $\B_n$ is the bridge attached to $\C$. On the other hand, if $\wh{\B}_n$ is a smaller bridge connecting the different sides of the bridge $\B_n$, let $\B_n\cup \wh{\B}_n$ be the handle $\h_n$. Then $\A\cup \h_n$ would be a cylinder with a handle (See Figure \ref{handle}).

\begin{figure}[b]
\begin{center}
$\begin{array}{c@{\hspace{.2in}}c}

\relabelbox  {\epsfxsize=2.5in \epsfbox{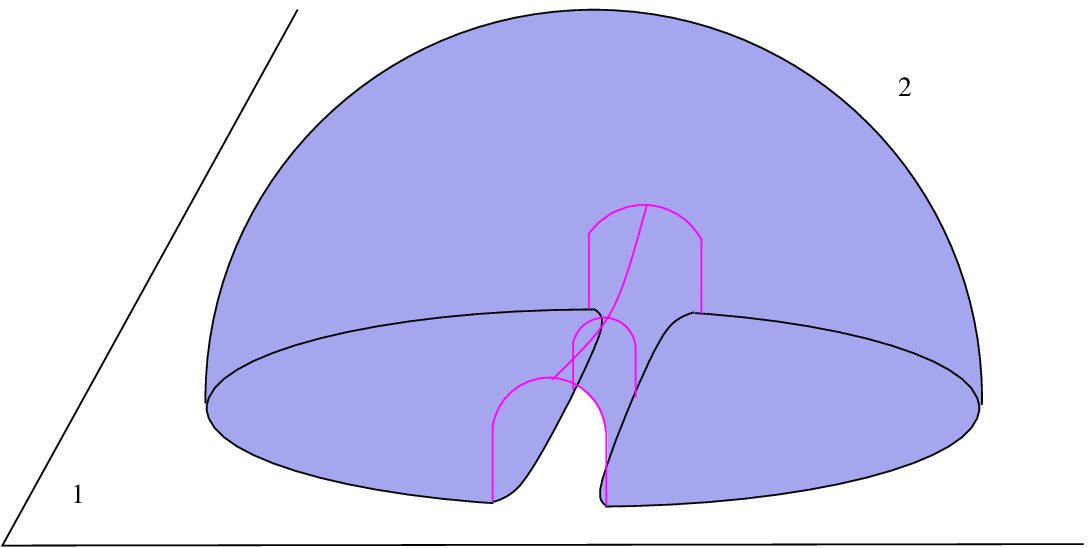}} \relabel{1}{$S^2_\infty$} \relabel{2}{$\Sigma_1\sharp \B_1$}  \endrelabelbox &

\relabelbox  {\epsfxsize=2.5in \epsfbox{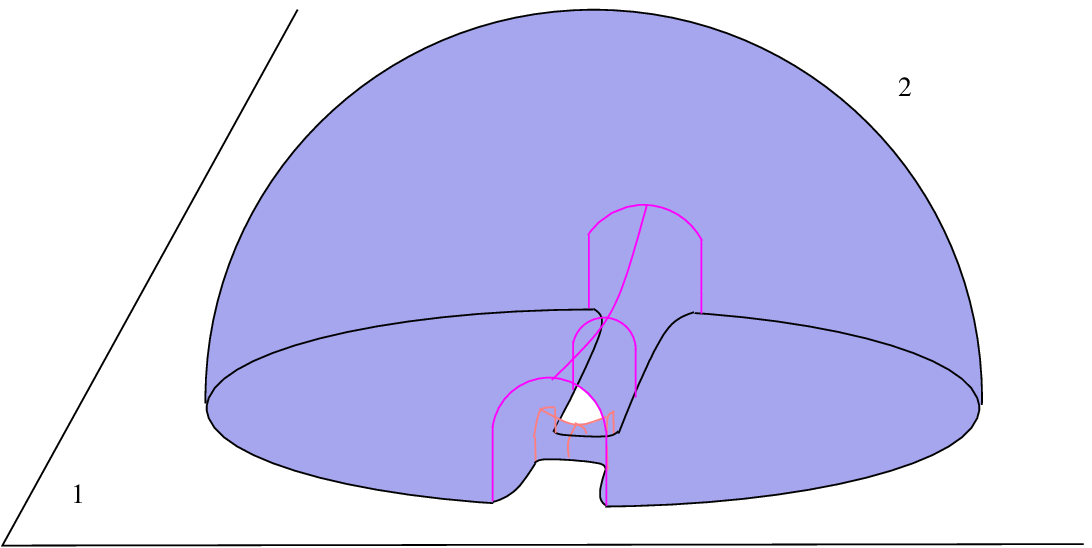}} \relabel{1}{$S^2_\infty$} \relabel{2}{$\Sigma_1\sharp\h_1$}    \endrelabelbox \\ [0.4cm]
\end{array}$

\end{center}
\caption{\label{caps} \small $\Sigma_1$ is a uniquely minimizing $H$-surface where $\PI\Sigma_1$ is a round circle. If $S_2-S_1$ contains a pair of pants, we attach one bridge $\B_1$ along $\beta_1$ to $\Sigma_1$, and get $\Sigma_2=\Sigma_1\sharp\B_1$ (left). If $S_2-S_1$ contains a cylinder with a handle, we attach two bridges successively to $\Sigma_1$ and get $\Sigma_2=\Sigma_1\sharp\h_1$ (right).  }
\end{figure}

Recall that if $\Sigma$ is an orientable surface of genus $g(\Sigma)$ with $k(\Sigma)$ boundary components, then its Euler Characteristic $\chi(\Sigma)=2-2g(\Sigma)-k(\Sigma)$.  Notice that by attaching a bridge $\B_n$, we increase the number of boundary components of $S_n$ by $1$ and decrease the euler characteristic by $1$, i.e. $\sharp(\partial S_{n+1})=\sharp(\partial S_n)+1$ and $\chi(S_{n+1})=\chi(S_n)-1$. Hence, $g(S_n)=g(S_{n+1})$ where $g(.)$ represents the genus of the surface. Similarly by attaching a handle $\h_n$ to $S_n$, we keep the number of boundary components same, but decrease the euler characteristic by $2$, i.e. $\sharp(\partial S_{n+1})=\sharp(\partial S_n)$ and $\chi(S_{n+1})=\chi(S_n)-2$. This implies $g(S_{n+1})=g(S_{n})+1$ with the same number of boundary components.

We start the construction with a minimizing $H$-plane $\Sigma_1$ (a spherical cap) in $\BHH$ bounding a round circle $\Gamma_1$ in $\Si$. Hence, $\overline{\Sigma}_1 \simeq S_1$. Now, we continue inductively (See Figure \ref{caps}). Assume that $S_{n+1}-S_n$ contains a pair of pants. Let the pair of pants attached to the component $\gamma$ in $\partial S_n$. Let $\gamma'$ be the corresponding component of $\Gamma_n =\PI\Sigma_n$. By construction, $\gamma'$ bounds a disk $D$ in $\Si$ with $D\cap\Gamma_n = \gamma'$. Let $\beta_n$ be a smooth arc segment in $D$ with $\beta_n\cap\Gamma_n = \partial\beta_n\subset \gamma'$, and $\beta_n \perp \gamma'$. Now, as $\Sigma_n$ is uniquely minimizing $H$-surface, and $\beta_n$ satisfies the conditions, by using the Theorem \ref{bridge}, we get a uniquely minimizing $H$-surface $\Sigma_{n+1}$ with $\overline{\Sigma}_{n+1}\simeq S_{n+1}$. Note also that by Theorem \ref{bridge}, we can choose the bridge along $\beta_n$ as thin as we want. Hence, in the Poincare ball model, we can get an increasing sequence $r_n\nearrow\infty$ such that $B_{r_n}(0)\cap \Sigma_{n+1} \simeq S_n$ and $B_{r_{n+1}}(0)\cap \Sigma_{n+1}\simeq S_{n+1}$.

Now, assume that $S_{n+1}-S_n$ contains a cylinder with a handle. Again, let $\gamma$ be the component of $\partial S_n$ where the cylinder with handle attached, and let $\gamma'\subset \Si$ be the corresponding component in $\PI \Sigma_n$. Let $D$ be the disk in $\Si$ with $\partial D=\gamma'$ and $D\cap\Gamma_n = \gamma'$. Like before, let $\beta_n$ be a smooth arc segment in $D$ with $\beta_n\cap\Gamma_n = \partial\beta_n\subset \gamma'$, and $\beta_n \perp \gamma'$. Now, by Theorem \ref{bridge}, we get a uniquely minimizing $H$-surface $\Sigma_{n+1}'$. Again, by choosing the bridge sufficiently thin, we can make sure that $B_{r_n}\cap \Sigma_{n+1}' \simeq S_n$. Now, let $\beta_n'$ be the small smooth arc in $D$ connecting the opposite sides of the bridge along $\beta_n$. Similarly, by using Theorem \ref{bridge}, we add another tiny bridge along $\beta_n'$ to $\Sigma_{n+1}'$ and get a uniquely minimizing $H$ surface $\Sigma_{n+1}$ where $\overline{\Sigma}_{n+1}\simeq S_{n+1}$. Like before, we can find sufficiently large $r_{n+1}>r_n$ with $B_{r_n}(0)\cap \Sigma_{n+1} \simeq S_n$ and $B_{r_{n+1}}(0)\cap \Sigma_{n+1}\simeq S_{n+1}$.

Hence, we get a sequence of uniquely minimizing $H$-surfaces $\Sigma_n$ in $\BHH$ such that for $r_n\nearrow \infty$, if $m>n$, $B_{r_n}(0)\cap \Sigma_m \simeq S_n$. By using the techniques in \cite{MW} and a diagonal sequence argument, we get a limiting surface $\Sigma$ in $\BHH$ where the convergence is smooth on compact sets.

$\Sigma$ is a minimizing $H$-surface in $\BHH$ as being limit of minimizing $H$-surfaces in $\BHH$. Moreover, $\Sigma\simeq S$ as $\Sigma\cap B_{r_n}(0) \simeq S_n$ as the convergence is smooth. Finally, $\Sigma$ is properly embedded in $\BHH$ as for any compact set $K\subset \BHH$, there exists $r_n>0$ with $K\subset B_{r_n}(0)$, and $B_{r_n}(0)\cap \Sigma\simeq S_n$ which is compact. The proof follows.
\end{pf}

\begin{rmk} Notice that to apply the bridge principle proved in the previous section, one needs that the original curve must bound a unique minimizing $H$-surface. Hence, in order to add the bridges successively, one needs to get a curve which bounds a unique minimizing $H$-surface after adding the bridge. This is the main idea here, and that is why we need the Theorem \ref{bridge} to give a uniquely minimizing $H$-surface after attaching the bridge.
\end{rmk}

\begin{rmk} If $S$ has infinite topology, then by following the arguments in \cite[Theorem 4.1]{MW}, it can be showed that the distinct ends of corresponding $H$-surface $\Sigma$ are disjoint. Similarly by following the arguments in \cite[Theorem 4.4]{MW}, it might be possible to show that $\PI \Sigma$ is a smooth curve except at one point. However, one might need more control in the asymptotic boundary of the surfaces, as our bridge principle completely changes the boundary curve at infinity unlike \cite{MW} where the boundary curve only changes near the bridge.
\end{rmk}

\section{Non-Properly Embedded Minimal Surfaces in $\BHH$}

In this section, we will show that any open orientable surface $S$ can also be nonproperly embedded in $\BHH$ as a minimal surface. The basic idea is by taking the  minimal surface $\Sigma_1$ in $\BHH$ with the desired topological type constructed in previous section, and the nonproperly embedded minimal plane $\Sigma_2$ in $\BHH$ constructed in \cite{Co3}, and "placing" a bridge between $\Sigma_1$ and $\Sigma_2$. For this construction, first we need a generalization of the bridge principle at infinity (Theorem \ref{bridge}) for area minimizing surfaces.

\subsection{A generalization of bridge principle at infinity for area minimizing surfaces}

Recall that Martin and White's bridge principle applies to uniquely minimizing surfaces in $\BHH$ \cite{MW}. In particular, if a collection of smooth simple closed curves $\Gamma$ in $\Si$ bounds a unique absolutely area minimizing surface $\Sigma$ (not necessarily connected) in $\BHH$, then for any smooth closed arc $\beta$ is $\Si$ with $\beta\cap \Gamma=\partial \beta$ where $\beta$ meets $\Gamma$ orthogonally, then there exists a unique minimizing surface $\wh{\Sigma}$ with $\wh{\Sigma}$ is "close" and homeomorphic to $\Sigma\cup N(\beta)$.

In this part, we will generalize this result, in particular Theorem \ref{bridge} for $H=0$, and this will be the key component of the construction of nonproperly embedded minimal surfaces in $\BHH$ with arbitrary topology.

\begin{defn} Let $\Sigma_1$ and $\Sigma_2$ be two complete uniquely minimizing surfaces in $\BHH$ with $\PI \Sigma_i = \Gamma_i$ where $\Gamma_i$ is a smooth collection of disjoint simple closed curves in $\Si$. If there exists a simple closed curve $\beta$ in $\Si$ such that $\Si-\beta=\Delta^+\cup \Delta^-$ and $\Gamma_1\subset \Delta^+$ and $\Gamma_2\subset\Delta^-$, then we will call $\Sigma_1$ and $\Sigma_2$ are {\em separated}.
\end{defn}

\begin{rmk} Notice that if $\Sigma_1$ and $\Sigma_2$ are separated, then $\Gamma_1$ and $\Gamma_2$ are disjoint, and hence, $\Sigma_1$ and $\Sigma_2$ are disjoint. On the reverse direction, unfortunately $\Gamma_1\cap\Gamma_2 =\emptyset$ and $\Sigma_1\cap\Sigma_2=\emptyset$ does not necessarily implies that $\Sigma_1$ and $\Sigma_2$ are separated, e.g. consider two area minimizing catenoids with the same rotation axis. Notice also that if $\Gamma_1$ and $\Gamma_2$ are disjoint simple closed curves (one component), then $\Sigma_1$ and $\Sigma_2$ are automatically separated.
\end{rmk}

Now assume that $\Sigma_1$ and $\Sigma_2$ are separated. Let $\alpha$ be a smooth closed arc in $\Si$ connecting $\Gamma_1$ and $\Gamma_2$ with $\alpha\cap (\Gamma_1\cup\Gamma_2)=\partial \alpha$ and $\alpha\bot \Gamma_i$. By using the notation at the beginning of Section 3, let $\wh{\Gamma}=\Gamma_1\cup\Gamma_2$ and define $Y$ as the one side of the neighborhood of $N(\wh{\Gamma}\cup\alpha)$, i.e. $Y=N^+(\wh{\Gamma})\cup N^+(\alpha)$. Let $\{\wh{\Gamma}_t \ | \ t\in(0,\delta) \}$ be the foliation of $Y$ where $\wh{\Gamma}_t\to \Gamma_1\cup\Gamma_2\cup\alpha$ as $t\searrow 0$.

\begin{thm} \label{genbridge} Let $\Sigma_1$ and $\Sigma_2$ be two uniquely minimizing surfaces in $\BHH$ with $\PI\Sigma_i=\Gamma_i$. Assume that $\Sigma_1$ and $\Sigma_2$ are separated. Let $\alpha$ be a smooth closed arc in $\Si$ connecting $\Gamma_1$ and $\Gamma_2$ with $\alpha\cap (\Gamma_1\cup\Gamma_2)=\partial \alpha$ and $\alpha\bot \Gamma_i$. Then for sufficiently small $t$, $\wh{\Gamma}_t$ bounds a minimal surface $\wh{\Sigma}_t$ in $\BHH$ which is homeomorphic to $\Sigma_1\cup\Sigma_2\cup N_\e^+(\alpha)$. In particular, $\wh{\Sigma}_t$ is uniquely minimizing surface in a mean convex subspace $X$ of $\BHH$.
\end{thm}

\begin{pf} There are two steps in the proof. A priori, it is not known whether $\Sigma_1\cup\Sigma_2$ is uniquely minimizing in $\BHH$. $\Sigma_1\cup\Sigma_2$ may not even be area minimizing in $\BHH$, see Remark \ref{genrmk}. In the first step, we will construct a mean convex domain $X$ in $\BHH$ which looks like a neighborhood of a plane is removed but still connected with a very thin solid cylinder. In the second step, we will show that $\Sigma_1\cup\Sigma_2$ is uniquely minimizing in $X$, and by using Theorem \ref{bridge}, we finish the proof.

\vspace{.2cm}

\noindent {\bf Step 1 - Construction of the Mean Convex Subspace $X$ in $\BHH$:}

\vspace{.2cm}

{\bf Step 1a - Construction of Igloos:} Let $\beta$ separate $\Gamma_1$ and $\Gamma_2$ in $\Si$, i.e. $\Si-\beta=\Delta^+\cup \Delta^-$ where $\Gamma_1\subset \Delta^+$ and $\Gamma_2\subset\Delta^-$. $\Sigma_1$ separates $\BHH$ into open components, and let $\Omega_1$ be the component in $\BHH-\Sigma_1$ with $\PI\Omega_1\supset \beta$. Similarly, define $\Omega_2$ as the component in $\BHH-\Sigma_2$ with $\PI\Omega_2\supset \beta$. Let $\Omega=\Omega_1\cap\Omega_2$. Then, $\overline{\Omega}$ would be a mean convex subspace of $\BHH$ with $\partial \overline{\Omega}\subset \Sigma_1\cup\Sigma_2$. Let $\Delta=\PI\overline{\Omega}\subset\Si$. Then, $\beta\subset \Delta$ and $\partial\Delta\subset \Gamma_1\cup\Gamma_2$.

Consider the handlebody $M=\overline{\Omega}\cup\Delta$. By construction $\partial M$ is connected, and it is a genus $g$ surface for $g=g(\Sigma_1)+g(\Sigma_2)$. We claim that $\beta\subset\partial M$ is nullhomotopic in $M$. Let $\mu_1, \tau_1, \mu_2, \tau_2, .. , \mu_g,\tau_g$ be the generators of the $\pi_1(\partial M)$ where $\mu_i$ curve is the meridian of genus $i$, where $\tau_i$ curve is the inner circle of genus $i$. Hence, each $\mu_i$ is trivial in $\pi_1(M)$. Notice that $\beta$ is a separating curve in $\partial M$, and it represents the trivial cycle in $H_1(\partial M)$. This means $\beta$ is just product of some commutators in $\pi_1(M)$, i.e. $\beta=[\mu_{i_1},\tau_{i_1}]*[\mu_{i_2},\tau_{i_2}]*..*[\mu_{i_k},\tau_{i_k}]$ where $[\mu,\tau]=\mu*\tau*\mu^{-1}*\tau^{-1}$. Since each commutator $[\mu_{i},\tau_{i}]$ is nullhomotopic in $M$, then $\beta=[\mu_{i_1},\tau_{i_1}]*[\mu_{i_2},\tau_{i_2}]*..*[\mu_{i_k},\tau_{i_k}]$ is nullhomotopic in $M$.

Now, let $\beta^+$ and $\beta^-$ are very close curves to $\beta$ in $\Si$ in the opposite sides, say $\beta^+\cup\beta^-=\partial N_\epsilon{\beta}$ where $N_\epsilon{\beta}$ is the $\epsilon$ neighborhood of $\beta$ in $\Si$ for sufficiently small $\epsilon>0$. As $\beta\pm$ nullhomotopic in $M$, by \cite{MY}, we can define a sequence of least area disks $D^\pm_i$ in $\Omega$ with $\partial D^\pm_i\to \beta^\pm \subset \Si$. Then, by \cite{A2}, \cite{Ga}, we obtain two least area planes $P^+$ and $P^-$ in $\Omega$ in the limit. Note that $P^\pm$ are just minimal planes in $\BHH$. By construction, $P^\pm\cap\Sigma_i=\emptyset$, and indeed $P^+$ (also $P^-$) separates $\BHH$ into two components where $\Sigma_1$ and $\Sigma_2$ belongs to different components.

By the definition of $\beta$, $\alpha\cap\beta\neq \emptyset$. By modifying $\beta$ if necessary, we can assume that $\alpha\cap\beta$ consists of just one point, say $x$. let $\eta=\alpha\cap N_\epsilon(\beta)$. Then $\eta$ is a short subarc in $\alpha$ between $\beta^+$ and $\beta^-$. Let $S_\delta(\eta)=N_\delta(\eta)\cap N_\epsilon(\beta)$ be a thin strip in $\Si$ along $\eta$. Let $\partial S_\delta(\eta)-(\beta^+\cup\beta^-)=\eta^+\cup\eta^-$. Define $\tau_\delta=(\beta^+\cup\beta^- -\partial S_\delta(\eta))\cup (\eta^+\cup\eta^-)$ which is a simple closed curve in $\Si$. Now, we will construct a minimal plane $\Pi$ in $\BHH$ such that $\PI \Pi = \tau_\delta$ and $\Pi \sim P^+\sharp_\eta P^-$ by using the techniques in \cite{Co3}.

Let $\gamma^+$ and $\gamma^-$ be two round circles in $N_\epsilon(\beta)$ in the opposite sides of $\eta$. By choosing $\gamma^+$ and $\gamma^-$ sufficiently close, we can make sure that there is a spherical catenoid $\mathcal{C}$ in $\BHH$ with $\PI \mathcal{C}=\gamma^+\cup\gamma^-$. Let $T$ (tunnel) be the "small" component of $\BHH-\mathcal{C}$ where $\PI T$ is union of two small disks in $\Si$ bounding $\gamma^+$ and $\gamma^-$. Let $\Omega'=\Omega-T$. Notice that $\Omega'$ is mean convex and $\tau$ is nullhomotopic in $\Omega'$ as $P^+\cup P^-\cup S_\delta(\eta)$ is a disk  in $\Omega$ with boundary $\tau$. Hence, like before, we can define a sequence of least area disks $D_i$ in $\Omega'$ with $\partial D_i\to \tau$, and in the limit, we get a least area plane $\Pi$ in $\Omega'$. By construction, $\Pi$ is "close" to $P^+\ \cup P^-\cup S_\delta(\eta)$ \cite{Co3}.

Now, $\Pi$ separates $\BHH$ into two components, and let $\mathcal{I}$ be the component of $\BHH-\Pi$ where $\PI \mathcal{I}$ does not contain $\eta$. In particular, $\mathcal{I}$ looks like an igloo, eskimo house, with a very tiny door, in the upper half space model (See Figure \ref{igloo}). Let $X= \BHH-\mathcal{I}$. Then, $X$ is a mean convex subspace of $\BHH$ with $\Sigma_1\cup\Sigma_2\subset X$. Notice that $X$ looks like a two large balls (inside and outside of the igloo) connected with a very thin solid cylinder, say {\em neck} (the tiny doorway of igloo). Notice also that while $\Sigma_1$ is {\em inside} of the igloo, $\Sigma_2$ is in the {\em outside} of the igloo.

\begin{figure}[t]

\relabelbox  {\epsfxsize=3.5in

\centerline{\epsfbox{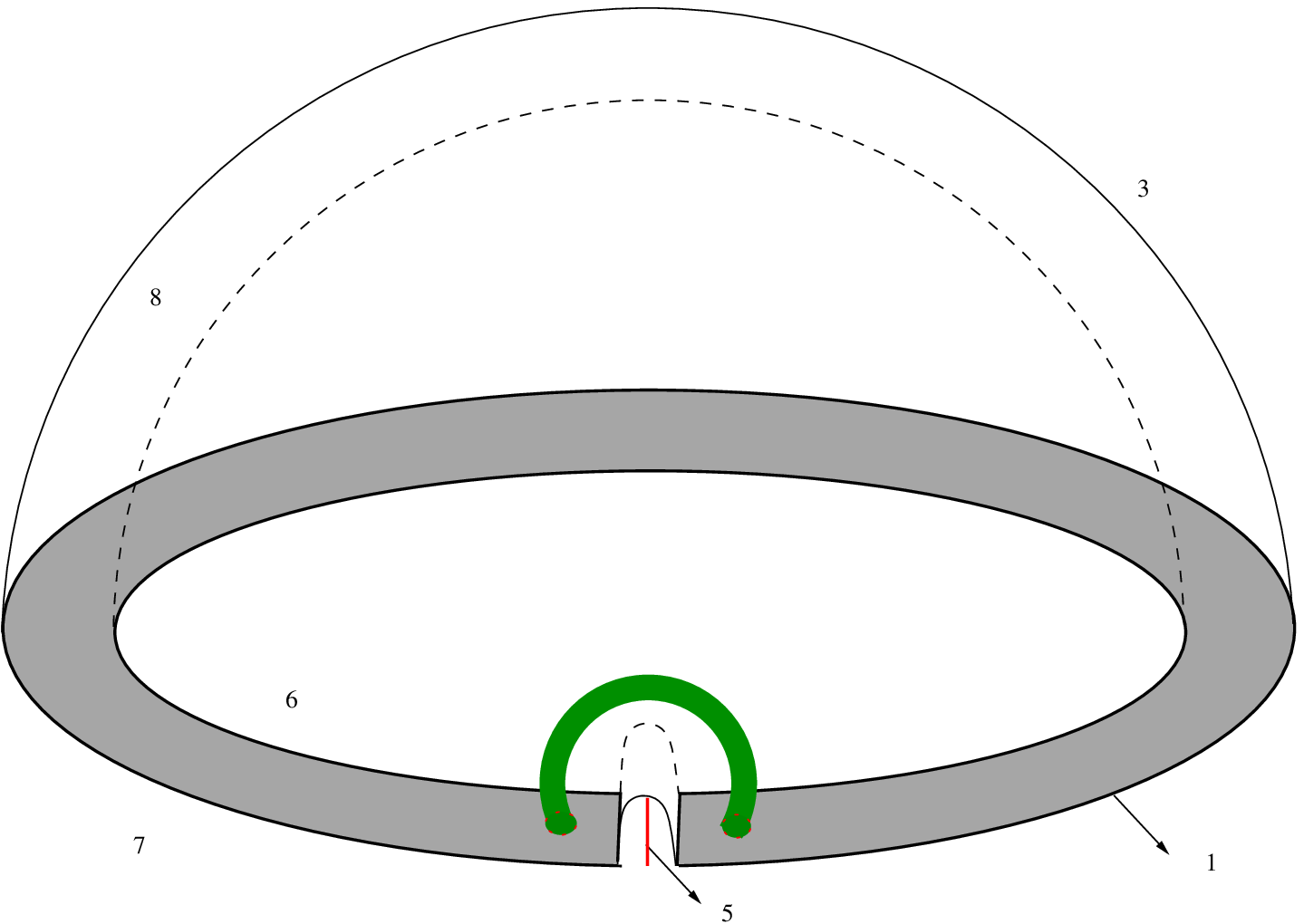}}}

\relabel{1}{$\tau$}

\relabel{3}{$\Pi$}

\relabel{5}{$\eta$}

\relabel{6}{$\beta^+$}

\relabel{7}{$\beta^-$}

\relabel{8}{$\mathcal{I}$}

\relabel{9}{$T$}

\endrelabelbox

\caption{\label{igloo} \small $\Pi$ is the least area plane in $X=\BHH-T$ where $\PI \Pi=\tau$. In particular, $\Pi=\p^+\sharp_{\eta}\p^-$ and  $\tau=\beta^+\sharp_\eta\beta^-$. }
\end{figure}

\vspace{.2cm}

{\bf Step 1b - Tiny Necks:} Notice that $\delta$ represents the width of the bridge along $\eta$ for $\Pi$, and as $\delta$ gets smaller, we get thinner necks in $X$. We will show that when the neck is sufficiently thin, there is no connected minimal surface going through the neck between the inside and outside of the igloo. In other words, {\em no minimal surface can pass through the neck}.

Now, we will use the upper half space model of $\BHH$. Let $q^+$ be the endpoint of $\eta$ in $\beta^+$, and let $q^-$ be the other endpoint of $\eta$ in $\beta^-$. Let $l^\pm$ be the tangent line of $\beta^\pm$ at $q^\pm$. Let $\mathcal{P}^\pm$ be the geodesic plane in $\BHH$ with $\PI \mathcal{P}^\pm = l^\pm$. Clearly, $\mathcal{P}^\pm$ cuts through the bridge near $\eta$ in $\Pi$. By choosing $\delta$ sufficiently small, and translating $\mathcal{P}^\pm$ into $\eta$ small amount, we can assume that $\mathcal{P}^\pm\cap\Pi$ contains a line $\lambda^\pm$ near neck such that one limit point of $\lambda^\pm$ is in $\eta^+$ and the other endpoint is in $\eta^-$, i.e. $\PI \lambda^\pm\subset \eta^+\cup\eta^-$.  Let $F^+$ be the component of $\mathcal{P}^+-\lambda^+$ near the bridge, i.e. $\PI F^+\cap \eta\neq \emptyset$. Similarly, let $F^-$ be the component of $\mathcal{P}^--\lambda^-$ near the bridge, i.e. $\PI F^-\cap \eta\neq \emptyset$. Clearly, $F^+$ and $F^-$ are area minimizing surfaces as they are subsurfaces of the geodesic planes $\mathcal{P}^+$ and $\mathcal{P}^-$ respectively. However, this does not automatically implies that their union $F^+\cup F^-$ is also an area minimizing surface. In this analogy, one might consider $F^+$ and $F^-$ are the inner and outer doors of the igloo.

For each $\delta=1/m$ (the thickness of the bridge), we can construct these planes, and say $F^\pm_m$ be the corresponding surfaces for $\delta=1/m$. Now, we will show that for sufficiently large $m$, $F^+_m\cup F^-_m$ is an area minimizing surface in $\BHH$.

Notice that if $F^+_m\cup F^-_m$ is not an area minimizing surface, then there are sufficiently large round circles $\zeta^+_m\subset F^+_m$ and $\zeta^-_m\subset F^-_m$ such that the area minimizing surface $\zeta^+_m\cup\zeta^-_m$ bounds in $\BHH$ is not the union of the disks $U^\pm_m \subset F^\pm_m$ with $\partial U^\pm_m = \zeta^\pm_m$, but a connected area minimizing surface $\A_m$ (e.g. annulus) with $\partial \A_m = \zeta^+_m\cup\zeta^-_m$. This is because if $U^+_m\cup U^-_m$ is not area minimizing, then any couple of larger disks $\wh{U}^+_m\cup\wh{U}^-_m$ with $U^\pm_m\subset\wh{U}^\pm_m$ is not area minimizing either. Hence, we  can choose $\zeta^+_m$ and $\zeta^-_m$ very large coaxial round circles in $F^+_m$ and $F^-_m$. Then, as $m\to \infty$, $d(\zeta^+_m,\zeta^-_m)\to \infty$. However, by \cite{Lo}, if $C_1$ and $C_2$ are distant circles in $\BHH$ with $d(C_1,C_2)>d_0$, then there is no connected minimal surface $S$ in $\BHH$ with $\partial S=C_1\cup C_2$. As $\A_m$ is a connected area minimizing -hence minimal- surface, this is a contradiction. This shows that for sufficiently large $m>0$, $F^+_m\cup F^-_m$ is an area minimizing surface in $\BHH$.

Hence, we fix a sufficiently large $m>0$ with $\delta=1/m$ for the mean convex subspace $X$ such that $F^+\cup F^-$ is area minimizimg.

\vspace{.2cm}

\noindent {\bf Step 2:} $\Sigma_1\cup\Sigma_2$ is uniquely minimizing surface in $X$.

\vspace{.2cm}

\noindent {\em Proof of Step 2:} Assume that there is an area minimizing surface $\Sigma'$ in $X$ different from $\Sigma_1\cup\Sigma_2$ with $\PI \Sigma' = \Gamma_1\cup\Gamma_2$. Since both $\Sigma_1$ and $\Sigma_2$ are uniquely minimizing in $\BHH$ by assumption, $\Sigma'$ must have a component $S$ such that $\PI S \nsubseteq \Gamma_1$ and $\PI S \nsubseteq \Gamma_2$. In other words, at least one end of $S$ is in $\Gamma_1$ and at least one end of $S$ is in $\Gamma_2$. This shows that $S$ must go through the neck region of $X$ near $\eta$.

Recall that $F^\pm=\mathcal{P}^\pm\cap X$, and $F^+\cup F^-$ is an area minimizing surface, i.e. any compact subsurface in $F^+\cup F^-$ is area minimizing. Since both $F^+$ and $F^-$ separates $X$, by construction $S\cap F^+\neq \emptyset$ and $S\cap F^-\neq \emptyset$. Since they are all area minimizing, the intersection must be a collection of closed curves, say $S\cap F^\pm=\sigma^\pm$. Let $S'$ be the compact subsurface of $S$ between $F^+$ and $F^-$, i.e. $\partial S' = \sigma^+\cup \sigma^-$. Let $D^+$ be the collection of disks in $F^+$ with $\partial D^+=\sigma^+$. Similarly, define $D^-$. As $S$ and $D^+\cup D^-$ are both area minimizing with same boundaries, $|S'|=|D^+|+|D^-|$ where $|.|$ represents the area.

Let $\wh{D}^+$ be a large disk in $F^+$ with $D^+\subset \wh{D}^+$. Similarly, define $\wh{D}^-\subset F^-$. Since $F^+\cup F^-$ is area minimizing surface, so is $\wh{D}^+\cup\wh{D}^-$. Define a new surface $\Sigma=(\wh{D}^+-D^+)\cup(\wh{D}--D^-)\cup S'$. Hence, $$|\Sigma|=(|\wh{D}^+|-|D^+|)+(|\wh{D}^-|-|D^-|)+|S'|=|\wh{D}^+|+|\wh{D}^-|$$

As $\wh{D}^+\cup\wh{D}^-$ is an area minimizing surface, and $\Sigma$ has the same area with the same boundary, $\Sigma$ is an area minimizing surface, too. However, $\Sigma$ has singularity along $\sigma^+\cup\sigma^-$. This contradicts to the regularity theorem for area minimizing surface \cite{Fe}. This shows that such an $S$ cannot exist, and the proof of Step 2 follows. \hfill \square

Finally, since $\Sigma_1\cup\Sigma_2$ is uniquely minimizing in the mean convex subspace $X$, by using Theorem \ref{bridge}, we obtain a uniquely minimizing surface $\wh{\Sigma}_t$ in $X$, which is homeomorphic to $\Sigma_1\cup\Sigma_2\cup N_\e^+(\alpha)$. The proof of the theorem follows.
\end{pf}

\begin{rmk} \label{MCbridge} Notice that in the proof, we can start with uniquely minimizing surfaces in a mean convex subspace $X$ of $\BHH$. Assume that $\Sigma_1$ and $\Sigma_2$ are uniquely minimizing surfaces in $X$, and they are separated in $X$ by a curve $\beta$ in $\PI X$ which is nullhomotopic in $X$. Further assume that the bridge $\alpha$ is in $\PI X$. Then the whole proof goes through, and we obtain a uniquely minimizing surface $\wh{\Sigma}=\Sigma_1\sharp_\alpha\Sigma_2$ in $X'=X-\mathcal{I}_\beta$ where $\mathcal{I}_\beta$ is the igloo over $\beta$ in $X$.
\end{rmk}

%\begin{rmk} An alternative way to prove the claim above is as follows. Let $S_m$ be the component of $\Pi_1^m$ separated by $l_1^m$ and $l_2^m$. Hence, $F_1^m\cup F_2^m\cup S_m$ separates a region $\Omega_m$ from $\BHH$. By rescaling $\Omega_m$ using the dilating isometries, $S_m$ limits to the infinite minimal strip, while $F_1^m$ and $F_2^m$ are pushed to infinity. Being far away from each other, $F_1^m\cup F_2^m$ must be area minimizing.\end{rmk}

\begin{rmk} \label{genrmk} This is an important generalization of Martin and White's bridge principle at infinity, as most of the time, the union of two uniquely minimizing surfaces in $\BHH$ may not be uniquely minimizing. Indeed, the union $\Sigma_1\cup\Sigma_2$ may not be area minimizing anymore, e.g. let $\Sigma_1$ and $\Sigma_2$ be two disjoint geodesic planes in $\BHH$ which are very close to each other. If they are sufficiently close, then the absolutely area minimizing surface for the union of their asymptotic boundary will not be the pair of geodesic planes, but instead it will be a spherical catenoid \cite{Wa}.
\end{rmk}

\begin{rmk} Notice that in the construction of $\Pi$, we used the least area planes $P^+$ and $P^-$ in $\Omega$, instead of the least area planes in $\BHH$. This is because the least area planes $P^\pm$ are disjoint from $\Sigma_1$ and $\Sigma_2$ by construction. However, the least area planes in $\BHH$ might intersect $\Sigma_1$ and $\Sigma_2$, which completely fails the construction. Hence, this choice is very important for the construction of the igloos, as it makes sure that the igloo $\mathcal{I}$ is disjoint from the surfaces $\Sigma_1$ and $\Sigma_2$, and $\Sigma_1\cup\Sigma_2\subset X=\BHH-\mathcal{I}$.
\end{rmk}

\newpage

\subsection{The Construction of Nonproperly Embedded Minimal Surfaces} \

\vspace{.2cm}

In this section we will construct nonproperly embedded minimal surfaces in $\BHH$ with arbitrary topology. In particular, we will show the following theorem:

\begin{thm} Any open orientable surface can be nonproperly embedded in $\BHH$ as a minimal surface. 
\end{thm}
\begin{pf} First, we give a short outline of the proof, and set the notation. Then, we proceed with the proof of the theorem.

\vspace{.2cm}

\noindent {\em Outline:} Let the open orientable surface $S$ be given. Let $\Sigma_1$ be the area minimizing surface in $\BHH$ which is homeomorphic to $S$ by Theorem \ref{main} and \cite{MW}. Let $\Sigma_2$ be the nonproperly embedded minimal plane in $\BHH$ by \cite{Co3}. Further assume that $\Sigma_1$ and $\Sigma_2$ are far away from each other, and $S^1_n\to \Sigma_1$ and $S^2_n\to \Sigma_2$ are the surfaces in the construction of $\Sigma_1$ and $\Sigma_2$. To construct nonproperly embedded minimal surface $\wh{\Sigma}$ with $\wh{\Sigma}\simeq S$, we will alternate the steps in these constructions, and define a new sequence $\{T_n\}$ of complete minimal surfaces, which is roughly $T_{2n} =S^1_n\sharp_\mu S^2_n$ where $\mu$ is the bridge between $\Sigma_1$ and $\Sigma_2$. Then, we show that $T_n\to \wh{\Sigma}$ is the minimal surface where $\wh{\Sigma}=\Sigma_1\sharp_\mu \Sigma_2$. Hence, $\wh{\Sigma}$ will have the same topological type with $\Sigma_1\simeq S$, and it will be nonproper because of $\Sigma_2$.

\vspace{.2cm}

\noindent {\em Notation and Setup:} Let $S$ be an open orientable surface. As in the previous section, let $S_1\subset S_2 \subset\ ....\ \subset S_n  \subset \ ...$ be a simple exhaustion of $S$ given by \cite{FMM}. Now, let $\wh{S}_n$ be the corresponding uniquely area minimizing embedding of $S_n$ into $\BHH$. Recall that in the construction in previous section, if $S_{n+1}-S_{n}$ contains a pair of pants, we are adding a suitable "bridge at infinity" to $\wh{S}_{n}$ in order to get $\wh{S}_{n+1}$, and if $S_{n+1}-S_{n}$ contains a cylinder with a handle, then we are adding "two bridges at infinity successively" to $\wh{S}_{n}$ in order to get $\wh{S}_{n+1}$. Without loss of generality, let $\beta_n$ represents this process dictated by the simple exhaustion, and say $\wh{S}_{n+1}=\wh{S}_{n}\sharp\beta_n$ for any $n$, i.e. $\beta_n$ represents a bridge if $S_{n+1}-S_{n}$ contains a pair of pants, and $\beta_n$ represents consecutive two bridges if $S_{n+1}-S_{n}$ contains a cylinder with handle (See Figure \ref{handle}).

To recall the construction of a nonproperly embedded plane $\Sigma_2$ in \cite{Co3}, let $\p_n$ be the geodesic plane where $\PI \p_n$ be the round circle $\gamma_n$ in $\Si$ with radius $1+1/n$ with center $(0,0,0)$ (upper half space model). Let $\p$ be the geodesic plane where $\PI \p$ is the round circle with radius 1 with center $(0,0,0)$. Clearly, $\p_n\to \p$. Now, we define minimal planes $E_n$ with $E_n=\p_1\sharp_{\alpha_1}\p_2\sharp_{\alpha_2}\ ... \ \sharp_{\alpha_{n-1}}\p_n$ where $\sharp_{\alpha_n}$ represents a bridge along $\alpha_n$ at infinity (See Figure \ref{planes}). However, the construction of these bridges is very different from the one in this paper.

\begin{figure}[t]
\begin{center}
$\begin{array}{l@{\hspace{.2in}}l}

\relabelbox  {\epsfxsize=2in \epsfbox{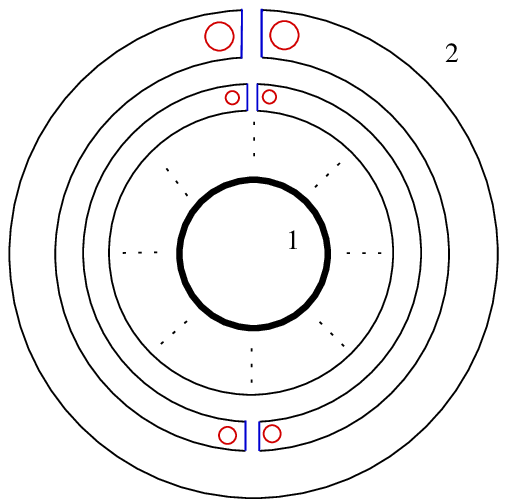}} \relabel{1}{$\gamma$} \relabel{2}{$\gamma_1$}  \endrelabelbox &

\relabelbox  {\epsfxsize=2.5in \epsfbox{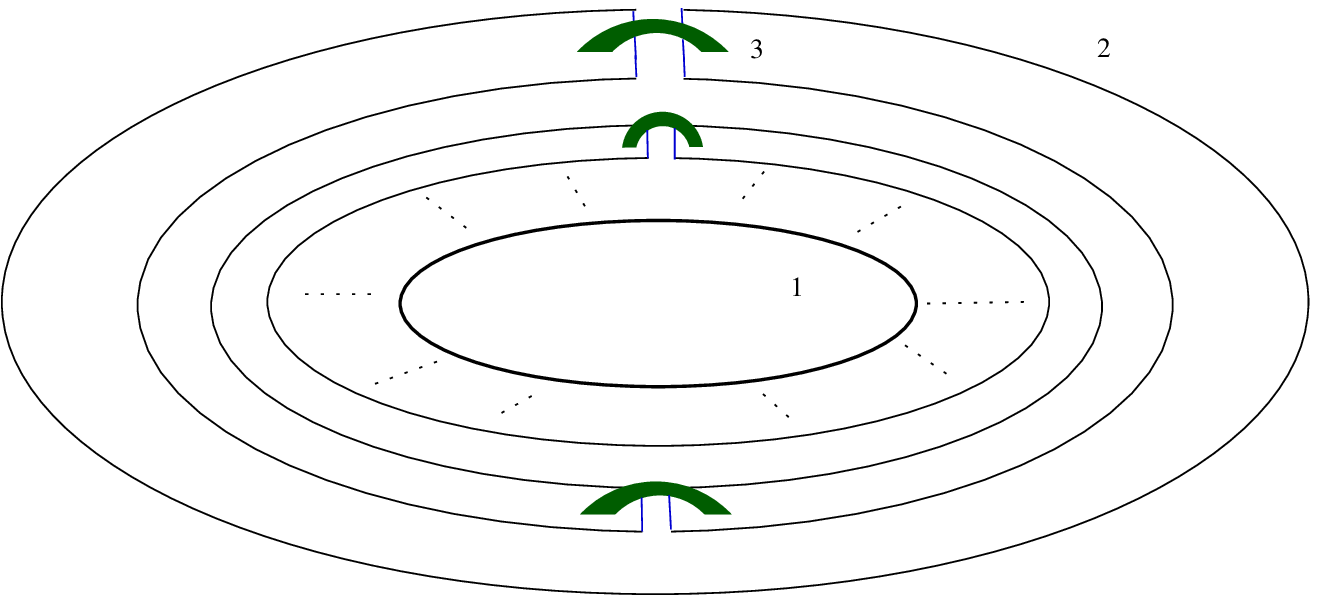}} \relabel{1}{$\gamma$} \relabel{2}{$\gamma_1$}   \endrelabelbox \\
\end{array}$

\end{center}
\caption{\label{planes} \small $\{\gamma_n\}$ is a sequence of round circles in $\Si$ where $\gamma_n\to\gamma$. $\alpha_n$ is the bridge connecting $\gamma_n$ and $\gamma_{n+1}$ (blue segments). The small red circles are the bases of the tunnels $\mathcal{T}_n$ (the green tubes in the right) which goes over the bridge $\alpha_n$. }
\end{figure}

Roughly, one needs to drill out a tunnel $\mathcal{T}_n$ which is the region inside a minimal catenoid in $\BHH$ where its ends are small circles in the opposite sides of the bridge $\alpha_n$. Then, $E_n$ is a least area plane in the mean convex subspace $X_n=X_{n-1}-\mathcal{T}_n$ with $\PI E_n=\Gamma_n$ where $\Gamma_n=\gamma_1\sharp\gamma_2\sharp...\sharp\gamma_n$. Notice that while $\Sigma_2$ is a least area plane in the mean convex subspace $X_\infty$ of $\BHH$, it is just a minimal plane in $\BHH$.

%Recall that to employ the bridge principle at infinity (Theorem \ref{bridge}), we need a uniquely minimizing surface to start with, while $E_n$ are not even area minimizing in $\BHH$. We will overcome this problem, by showing that they are area minimizing in a smaller subspace $X_n$ of $\BHH$, and get uniquely minimizing surfaces in $X_n$ for the bridge principle at infinity.

\vspace{.2cm}

\noindent {\em Construction of the Sequence $T_n$:} Now, we are inductively building the sequence of minimal surfaces $T_n$ in $\BHH$, which will give us the desired nonproperly embedded minimal surface $\wh{\Sigma}$, i.e. $T_n\to\wh{\Sigma}$. Note that we will construct the first four surfaces $T_1, T_2, T_3, T_4$ of the sequence explicitly. The construction of the remaining surfaces $T_n$ in the sequence will be clear.

In the construction above, we translate "right" the construction of $\Sigma_1$ by the parabolic isometry $\varphi_2(x,y,z)= (x+2,y,z)$. Hence,  $\wh{S}_1$ is the geodesic plane where $\PI \wh{S}_1$ is the circle $\eta$ of radius $1$ with center $(2,0,0)$. Similarly, we translate "left" the construction of $\Sigma_2$ by the isometry $\varphi_{-2}(x,y,z)= (x-2,y,z)$. Hence, $\p_1$ is the geodesic plane where $\PI \p_1$ is the circle $\gamma_1$ of radius $2$ with center $(-2,0,0)$.

Let $T_1=\wh{S}_1$. Let $\mu$ be the arc $[0,1]\times \{0\}\times\{0\}$ in $\Si$, connecting $\PI \wh{S}_1$ and $\PI \p_1$. Since both $\wh{S}_1$ and $\p_1$ are uniquely minimizing surfaces in $\BHH$, and they are separated (say by the round circle $\lambda_1$ with center $(2,0,0)$ of radius 3/2), we can use Theorem \ref{genbridge} to get a minimal surface $T_2=\wh{S}_1\sharp_\mu\p_1$. By Theorem \ref{genbridge}, note also that $T_2$ is a uniquely minimizing surface in the mean convex subspace $X_1$ of $\BHH$ with $X_1=\BHH-\mathcal{I}_1$ where $\mathcal{I}_1$ is the igloo over $\lambda_1$ (See Figure \ref{igloo}).

As $T_2$ is a uniquely minimizing surface in a mean convex subspace $X_1$, by applying Theorem \ref{bridge}, we get $T_3=T_2\sharp \beta_1$ where $\beta_1$ represents the collection of bridges, or handles (successive two bridges) in the construction of $\Sigma_1$ as in the proof of Theorem \ref{main}. Hence, $T_3$ is homeomorphic to $S_2$ in the simple exhaustion, and uniquely minimizing surface in $X_1$.

Now, we define $T_4$. Let $\alpha_1$ be the arc connecting $\gamma_1$ and $\gamma_2$ in the construction of $\Sigma_2$ (See Figure \ref{planes}). $T_3$ and $\p_2$ are uniquely minimizing surfaces in the mean convex subspace $X_0$, and they are separated by the round circle $\lambda_2$ with center $(-2,0,0)$ of radius $7/4$). Hence, we can apply the generalized version of the bridge principle at infinity (Theorem \ref{genbridge} and Remark \ref{MCbridge}) to $T_3\cup \p_2$ once again, and we get a uniquely minimizing surface $T_4=T_3\sharp_{\alpha_1}\p_2$ in a mean convex subspace $X_2=X_1-\mathcal{I}_2$. Here, $\mathcal{I}_2$ represents the igloo over $\lambda_2$.

After defining the first $4$ surfaces in the sequence, we can construct the remaining surfaces in the sequence inductively as follows.

By induction, $T_{2n}$ is uniquely minimizing in $X_n$. Hence, by applying Theorem \ref{bridge}, define $T_{2n+1}=T_{2n}\sharp\beta_n$ where $\beta_n$ represents the corresponding bridges, or handles in the construction of $\Sigma_1$. Hence, $T_{2n+1}$ is homeomorphic to $S_{n+1}$ in the simple exhaustion, and uniquely minimizing in $X_n$.

Define $T_{2n}$ as follows. By induction, $T_{2n-1}$ is uniquely minimizing in the mean convex subspace $X_{n-1}$. Since $\p_n$ is uniquely minimizing in $\BHH$, it is automatically minimizing in $X_{n-1}$. Notice that by convex hull property for any $n,m>0$ $\mathcal{I}_n\cap \p_m=\emptyset$. Hence, $T_{2n-1}$ and $\p_n$ are uniquely minimizing in $X_{n-1}$ and they are separated by the round circle $\lambda_n$ with center $(-2,0,0)$ of radius $1+\frac{2n+1}{2n(n+1)}$. Let $\alpha_{n-1}$ be the arc connecting $\gamma_{n-1}$ and $\gamma_n$ in the construction of $\Sigma_2$ (See Figure \ref{planes}). Then, by applying Theorem \ref{genbridge}, we obtain $T_{2n}=T_{2n-1}\sharp_{\alpha_{n-1}}\p_n$ which is a uniquely minimizing surface in $X_n=X_{n-1}-\mathcal{I}_n$ where $\mathcal{I}_n$ is the igloo over $\lambda_n$.

\vspace{.2cm}

\noindent {\em Nonproperly Embedded Minimal Surfaces with Arbitrary Topology:} Let $X_\infty = \bigcap_{n=1}^\infty X_n$ be the mean convex region in $\BHH$. As $\mathcal{I}_m\cap T_n=\emptyset$ for any $m>n$ by convex hull property, $T_n\subset X_\infty$. Since $X_\infty\subset X_n$ for any $n$, $T_n$ is a uniquely minimizing surface in $X_\infty$. Then the limit surface $\wh{\Sigma}=T_\infty$ is an area minimizing surface in $X_\infty$ and hence a minimal surface in $\BHH$. Clearly, $\wh{\Sigma}$ has the same topological type with the given surface $S$ by the construction, i.e. $\wh{\Sigma}\simeq \Sigma_1\sharp_\mu \Sigma_2$. $\wh{\Sigma}$ is nonproper as the closure of $\wh{\Sigma}$ is $\wh{\Sigma} \cup \p$ where $\p$ is the geodesic plane in $\BHH$ with $\PI \p$ is a round circle of radius $1$ and center $(-2,0,0)$. The proof follows.
\end{pf}

\section{Final Remarks}

In this paper, we first generalize Martin and White's result on the existence of complete area minimizing surfaces in $\BHH$ of arbitrary topological type. In particular, they showed that if $S$ is an open orientable surface, then there exists a complete proper embedding of $S$ into $\BHH$ as an area minimizing surface \cite{MW}. We generalize this result by showing that there exists a complete proper embedding of $S$ into $\BHH$ as an $H$-minimizing surface for $0\leq H<1$. Note that here $H=0$ corresponds to the area minimizing case.

When generalizing their result, our approach is mainly similar, but techniques are very different in some particular steps. In both papers, when constructing the topology of the given surface, the main tool is the bridge principle. In order to use this bridge principle, both approach needs the original surface to be uniquely minimizing to start with. Also, to apply this bridge principle again, the resulting surface after the bridge attached should be uniquely minimizing, too.

In order to ensure the uniqueness after the bridges attached, while Martin and White use the analytic tools, namely {\em $L^\infty$ stability} condition on the surfaces, we use the generic uniqueness tools developed in \cite{Co2}, which are more topological. On the other hand, in order to prove the resulting surface after the bridge attached has the desired topology, Martin and White uses strips and skillets idea from the the original bridge principle theory developed by White \cite{Wh}. In particular, they used these tools to show that there is no genus developed in the bridge when attaching. Similarly, in this paper, we followed their methods for the same step, and generalized their minimal strips, and skillets idea as $H$-strips and skillets in the appendix.

While in section 3, we showed the existence of properly embedded $H$-surfaces in $\BHH$ of arbitrary topological type, in the following section, we generalize Martin and White's result in a different direction. Especially after Colding and Minicozzi's proof of the Calabi-Yau Conjecture \cite{CM}, the nonproper embeddings of minimal surfaces became very interesting. We show that if $S$ is an open orientable surface, then there exists a complete {\bf nonproper} embedding of $S$ into $\BHH$ as a minimal surface. We show this by "placing a bridge" between the area minimizing surface of topological type of $S$ like above, and a minimal plane constructed in \cite{Co3}. First of all, unfortunately this surface is not area minimizing but just minimal in $\BHH$ by construction. It would be an interesting question whether there exists a nonproperly embedded area minimizing surface in $\BHH$ of arbitrary topological type.

On the other hand, while we can construct properly embedded $H$-surfaces of arbitrary topological type in $\BHH$, the same techniques do not apply to construct {\bf nonproperly embedded $H$-surfaces} in $\BHH$. In particular, in the construction above, we have this nonproperly embedded minimal plane, and we are attaching it via a bridge to the area minimizing surface of topological type of $S$. However, in $0<H<1$, a similar nonproperly embedded $H$-plane does not exist to start with. This is simply because the construction in \cite{Co3} does not apply to $0<H<1$ case, because of the orientation issues. In Section 4, a summary of this construction is given, and when we can attach the minimal planes $\p_n$ and $\p_{n+1}$ via bridge and get another minimal plane. However, for $H$-planes this is not possible. When we attach corresponding $\p^H_n$ and $\p^H_{n+1}$, the bridge does not connect the convex sides. In particular, when one end connects to an $H$-surface, the other end connects to $-H$-surface, hence the construction fails very seriously. On the other hand, Meeks, Tinaglia and the author showed the existence of the nonproperly embedded $H$-plane in $\BHH$ for $0\leq H<1$, which is an infinite strip spiraling between two $H$-catenoids \cite{CMT}. It might be possible to apply the construction above with this nonproperly embedded $H$-plane, which would show the existence of nonproperly embedded $H$-surfaces in $\BHH$ of arbitrary topological type.

One other very interesting question coming out of the construction of nonproperly embedded minimal surfaces is the a general {\em bridge principle at infinity  for complete, stable minimal surfaces} in $\BHH$. The bridge principle at infinity developed in \cite{MW}, or in this paper is just for  uniquely minimizing surfaces. One suspects that a more general version might be true. In particular, it is a very interesting question whether the bridge principle at infinity is true for globally stable minimal surfaces (or $H$-surfaces) in $\BHH$, i.e. if $\Sigma_1$ and $\Sigma_2$ are globally stable minimal surfaces in $\BHH$ with $\PI\Sigma_i = \Gamma_i$, and $\alpha$ is an arc in $\Si$ between $\Gamma_1$ and $\Gamma_2$, then is there a complete stable minimal surface $\wh{\Sigma}=\Sigma_1\sharp_\alpha\Sigma_2$ with $\wh{\Sigma}\sim \Sigma_1\cup\alpha\cup\Sigma_2$? It is reasonable to expect to use the tools (like igloo trick) in the proof of Theorem \ref{genbridge} to employ the techniques in {\em the original bridge principle for stable minimal surfaces} \cite{Wh}. Recall that to prove the original bridge principle for stable minimal surfaces $\Sigma_1$ and $\Sigma_2$ where $\partial \Sigma_i =\Gamma_i$ and $\alpha$ is an arc connecting $\Gamma_1$ and $\Gamma_2$, one first constructs a small mean convex neighborhood $N$ of $\Sigma_1\cup\alpha\cup\Sigma_2$ in the ambient space. Then, the area minimizing surface in $N$ bounding $\Gamma=\Gamma_1\sharp_\alpha\Gamma_2\subset \partial N$ is a minimal surface very close to $\Sigma_1\cup\alpha\cup\Sigma_2$ because of the choice of $N$. Hence, if one can construct the appropriate mean convex neighborhood $X$ of $\Sigma_1\cup\Sigma_2\cup\alpha$ in $\BHH$, and solve the Plateau problem in $X$ for $\Gamma=\Gamma_1\sharp_\alpha\Gamma_2\subset \PI X$, it would give the desired surface, and prove the bridge principle at infinity in full generality.

\section{Appendix: H-strips and H-skillets}

In this part, we will show that there is no genus developed in the bridge near infinity in Theorem \ref{bridge}. We will use the notation of Section 3. In particular, let $\Gamma_0$ be a collection of simple closed curves in $\Si$ which bounds a unique $H$-surface $T$. Let $\alpha$ be the bridge and $\Gamma_t$ be a foliation of positive part of the neighborhood $N_\epsilon(\Gamma\cup \alpha)$ with $\Gamma_t\to\Gamma_0\cup\alpha$ as $t\searrow 0$. Let $S_t$ be minimizing $H$-surface in $\BHH$ with $\PI S_t = \Gamma_t$. Then, as in section $3$, there exists a sequence $t_n\searrow 0$ with $S_{t_n} \to T$. Say $S_n=S_{t_n}$ and $\Gamma_n=\Gamma_{t_n}$. Let $\R_a=\{0\leq z \leq a\}$ in $\overline{\BHH}$. In this section, we will prove the following lemma:

\begin{lem} \label{nogenus} There exists $a>0$ and $N>0$ such that for any $n>N$, $S_n\cap R_{a}$ has no genus, i.e. $S_n\cap R_{a}\simeq \Gamma_0\times (0,a]$
\end{lem}

\begin{pf} Assuming that $H$-strips and $H$-skillets are uniquely minimizing $H$-surfaces (proved below), the proof is as follows. Similar to \cite{MW}, assume on the contrary that for any $a>0$, there exists a subsequence $S_n\cap R_{a}$ has genus. Then, let $\Delta_n$ be the component of $\BHH-S_n$ which contains the bridge $\alpha$. Since $S_n\cap R_a$ has genus, then $\Delta_n\cap R_a$ must be a nontrivial handlebody, i.e. it is not a $3$-ball. Hence, there must be a point $p_n$ in $S_n\cap R_{a}$ where the normal vector $v_{p_n}=<0,0,1>$ pointing inside $\Delta_n$.

%Similarly, there exists a point $q_n$ in $S_n\cap R_{a}$ where the normal vector $v_{q_n}=<0,0,-1>$ pointing inside $\Delta_n$.

Let $p_n=(x_n,y_n,z_n)$. Consider the isometry $\psi_n(x,y,z)=\frac{1}{z_n}(x-x_n,y-y_n,z)$ which is a translation by $-(x_n,y_n,0)$ first, and homothety by $\frac{1}{z_n}$ later. Then, consider the sequence of minimizing $H$-surfaces $S_n'=\psi_n(S_n)$ and $p_n'=\psi_n(p_n)=(0,0,1)$. Let $\Gamma_n'=\psi_n(\Gamma_n)=\PI S_n'$. After passing to a subsequence, we get the limits $S_n'\to S'$, $p_n'\to p'=(0,0,1)\in S'$, and $\Gamma_n'\to\Gamma'$. Note also that by construction the normal vector to $S'$ at $p'$ is $v_{p_n}\to v_p'=<0,0,1>$ pointing inside $\Delta'$.

Then like \cite{MW}, there are 4 possibilities. $\Gamma'$ is either a line, a T-shape, the union of two parallel lines or the boundary of a skillet. If $\Gamma'$ is a line or a T-shape, then $S'$ would be a half plane which makes $\theta_H$ angle with the $xy$-plane, i.e. $\Si$. Hence, the normal vector cannot be $<0,0,1>$ for any point in $S'$.

If $\Gamma'$ is the union of two straight lines, then $S'$ must be an $H$-strip for $-1<H<1$ as $H$-strips are uniquely $H$-minimizing. However, there is no normal vector $<0,0,1>$ on $H$-strips pointing inside $\Delta'$. Similarly, if $\Gamma'$ is the boundary of a skillet, then $S'$ would be an $H$-skillet for $-1<H<1$ as $H$-skillets are uniquely $H$-minimizing. Again, as there is no normal vector $<0,0,1>$ on $H$-skillets pointing inside $\Delta'$, this is a contradiction.
\end{pf}

Now, following \cite{MW} and \cite{Wh}, we define $H$-strips and $H$-skillets in $\BHH$, and show that $H$-strips and $H$-skillets are uniquely minimizing $H$-surfaces with special asymptotic boundaries.

\begin{figure}[b]
\begin{center}
$\begin{array}{c@{\hspace{.2in}}c}

\relabelbox  {\epsfxsize=2.5in \epsfbox{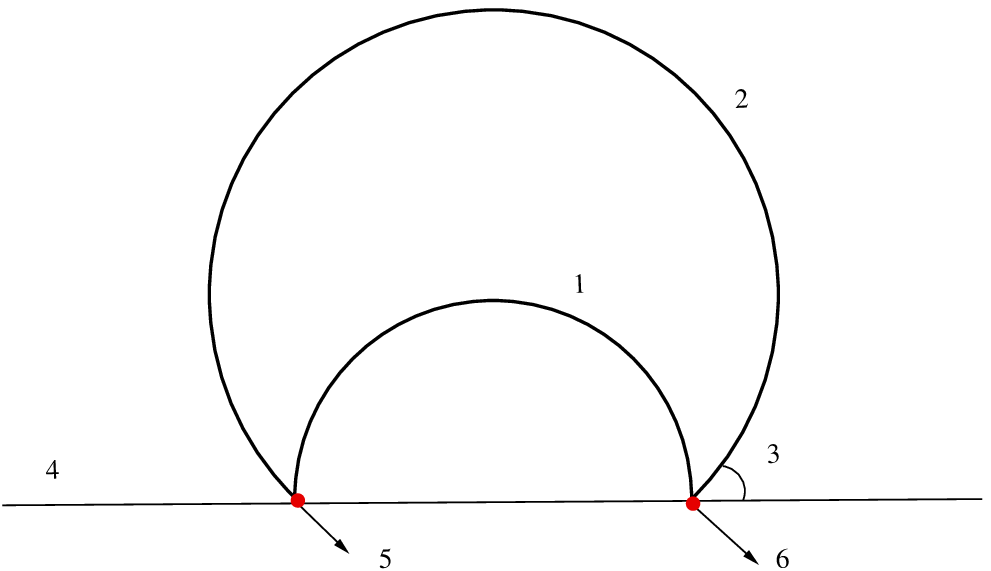}} \relabel{1}{$\eta_0$} \relabel{2}{$\eta_H$} \relabel{3}{$\theta_H$}  \relabel{4}{$\Si$} \relabel{5}{$\beta^-_\e$} \relabel{6}{$\beta^+_\e$}  \endrelabelbox &

\relabelbox  {\epsfysize=1.7in \epsfbox{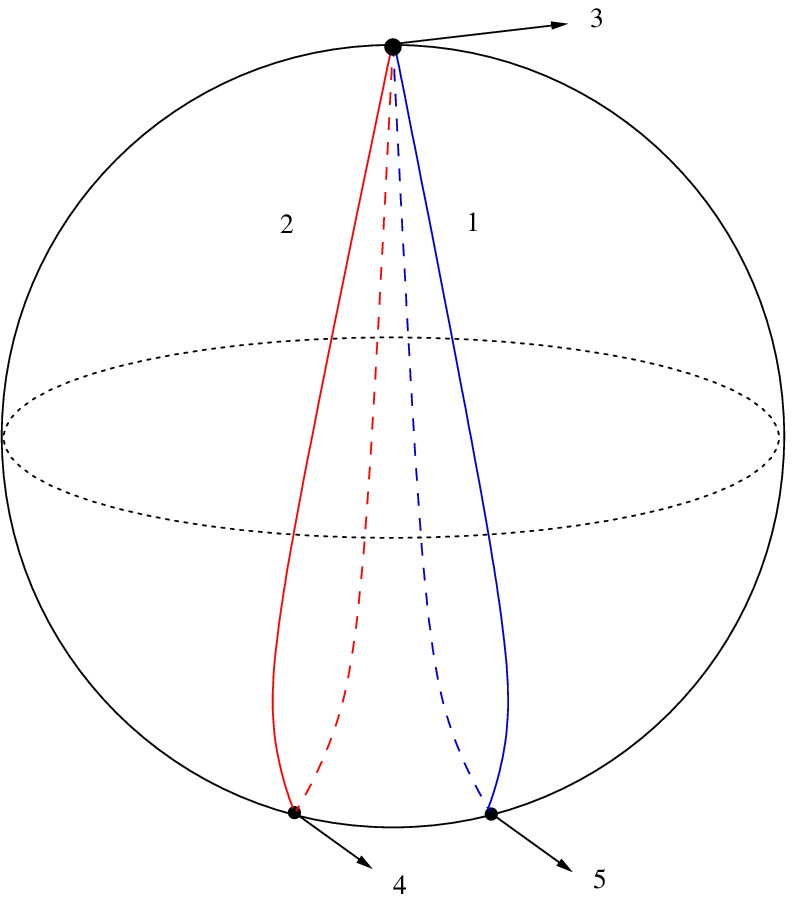}} \relabel{1}{$\beta^+_\e$} \relabel{2}{$\beta^-_\e$}  \relabel{3}{$\infty$} \relabel{4}{$p^-$} \relabel{5}{$p^+$}  \endrelabelbox \\ [0.4cm]
\end{array}$

\end{center}
\caption{\label{H-strip} \small In the left, if $\Sigma_H$ is an $H$-strip with $\PI \Sigma_H = \beta^+_\e\cup \beta^-_\e$ for $-1<H<1$, then $\Sigma_H=\eta_H\times \BR$ where $\eta_H$ is a smooth arc in $yz$-plane with endpoints $(\epsilon, 0)$ and $(-\epsilon,0)$, and which makes angle $\theta_H$ with $\Si$. In the right, the lines $\beta^+_\e$ and $\beta^-_\e$ are pictured in the Poincare ball model, where $p^\pm$ corresponds to points $(0,\pm\e ,0)$. }
\end{figure}

First, we define $H$-strips. We use the upper half space model for $\BHH$. Hence, $\Si=\{z=0\}\cup\{\infty\}$. With this notation, let $\beta_\epsilon$ be union of two straight lines parallel to $x$-axis in $xy$-plane in upper half space model, i.e. $\beta_\epsilon = \beta^+_\epsilon \cup \beta^-_\epsilon= \{ (x,\epsilon, 0)\} \cup \{(x,-\epsilon,0)\}$. Let $\Omega$ be the region in $\Si$ between these two lines, i.e. $\Omega=\{(x,y,0)\ | \ |y|\leq \epsilon\}$. Notice that $\beta_\epsilon$ is union of two round circles $\beta^\pm_\epsilon$ in $\Si$ where they touch each other at one point ($\infty$) in Poincare ball model (See Figure \ref{H-strip}). Let $\Sigma_\epsilon$ be a minimizing $H$-surface with $\PI \Sigma_\epsilon = \beta_\epsilon$ (Lemma \ref{existence}). We call $\Sigma_\epsilon$ an $H$-strip.\\

\noindent {\bf Claim:} $H$-strips are uniquely $H$-minimizing.

\vspace{.2cm}

\begin{pf} Fix $\epsilon$, and let $\beta_\epsilon = \beta$. By \cite{To}, there exists a minimizing $H$-surface $\Sigma$ with $\PI \Sigma =\beta$. First, notice that for sufficiently small $\epsilon$, $\Sigma$ is connected. To see that, if $\Sigma$ is not connected, then $\Sigma=\mathcal{P}_H^+\cup\mathcal{P}_H^-$ where $\mathcal{P}_H^\pm$ is the minimizing $H$-surface with $\PI \mathcal{P}_H^\pm=\beta^\pm$. Notice that as $\beta^+$ and $\beta^-$ are both round circles in $\Si$ (Poincare ball model), $\Sigma^+$ and $\Sigma^-$ are both spherical caps corresponding to $0\leq H<1$. In particular, in upper half space model, $\mathcal{P}_H^\pm$ would be half planes with $\mathcal{P}_H^+ = \{ y=\frac{H}{\sqrt{1-H^2}} z + \epsilon\}$ and $\mathcal{P}_H^- = \{ y=-\frac{H}{\sqrt{1-H^2}} z - \epsilon\}$ which makes angle $\theta_H$ with $\Si$ \cite{To}.

Now, let $\mathcal{C}_H$ be the spherical $H$-catenoid (see \cite{Go}) with $\PI \mathcal{C}_H = \gamma^+\cup \gamma^-$ where $\gamma^\pm$ are two small round circles in opposite sides of $\beta^+\cup \beta^-$, e.g. $\gamma^\pm$ is the circle of radius $0<\delta<\e$ with center $(0,\pm2\epsilon,0)$. In particular, $\mathcal{C}_0$ intersects $\mathcal{P}_0^\pm$ in two round circles $\tau^+\cup\tau^-$ which bounds two disks $D^+\cup D^-$ in $\mathcal{P}_0^+\cup\mathcal{P}_0^-$, and an annulus $A$ in $\mathcal{C}_0$. For sufficiently small $\epsilon>0$, and appropriate choice of $\delta$, $\mathcal{C}_0$ is a least area catenoid by \cite{Wa}. Hence, $|D^+|+|D^-| \geq |A|$ and it is easy to show that $\mathcal{P}^+_0 \cup \mathcal{P}^-_0$ is not area minimizing by a swaping argument. Similar comparison argument with the spherical $H$-catenoid $\C_H$ shows that $\mathcal{P}^+_H \cup \mathcal{P}^-_H$ is not a minimizing $H$-surface. Hence, the minimizing $H$-surface $\Sigma$ with $\PI \Sigma = \beta^+\cup\beta^-$ must be connected.

Now, assume that $\Sigma$ is not uniquely minimizing $H$-surface. Then by Lemma \ref{canonical}, there exists canonical minimizing $H$-surfaces $\Sigma^+$ and $\Sigma^-$ with $\PI \Sigma^\pm =\beta^+\cup \beta^-$. In particular, let $\Omega_i^-\subset \Si$ be an exhaustion of $\Omega$ by compact connected regions, i.e. $\Omega_1^-\subset \Omega_2^- \subset .. \Omega_i^-\subset ..$ with $\Omega =\bigcup_{i=1}^\infty \Omega_i$. Let $\partial \Omega_i^- = \alpha_i^-$, and $\Sigma_i^-$ be the minimizing $H$-surface with $\PI \Sigma_i^- = \alpha_i$ by Lemma \ref{existence}. Then, as $\alpha_i^- \to \beta^+\cup\beta^-$ by construction, there is a convergent subsequence $\Sigma_i^- \to \Sigma^-$ where $\Sigma^-$ is a minimizing $H$-surface with $\PI \Sigma^-=\beta^+\cup\beta^-$. Similarly, one can define $\Sigma^+$ by using a decreasing sequence of regions $\Omega^+_i$ in $\Si$ with $\Omega^+_{i+1}\supset\Omega^+_i\supset\Omega$ and $\Omega=\bigcap_{i=1}^\infty \Omega^+_i$. Define $\alpha_i^+$ and $\Sigma_i^+\to\Sigma^+$ similarly. Moreover, by Lemma \ref{canonical}, $\Sigma^+$ and $\Sigma^-$ are canonical and independent of the choices of $\{\Omega_i\}$.

Now, consider the parabolic isometry $\varphi_t$ of $\BHH$ which is a translation along $x$-axis, i.e. $\varphi_t(x,y,z) =(x+t, y, z)$. Clearly, $\varphi_t$ fixes $\beta^+$ and $\beta^-$ for any $t$, i.e. $\varphi_t(\beta^\pm)=\beta^\pm$. Let $\varphi_t(\Sigma^-)=\Sigma^-_t$. Clearly, $\PI\Sigma^-_t=\PI \Sigma^-$ for any $t$. On the other hand, $\Sigma^-_t$ is the limit of $\Sigma^t_i$ with $\PI\Sigma^t_i =\alpha^t_i=\partial \Omega^t_i=\partial \varphi_t(\Omega_i)$. However, $\Sigma^-$ is canonical, and it is independent of $\{\Omega_i\}$. Hence for any $t$, $\Sigma_t^-=\Sigma^-$. As $\Sigma^-$ is invariant under $\varphi_t$, this shows that $\Sigma^- = \eta^- \times\BR$ where $\eta^-$ is a smooth simple arc in $yz$-plane with endpoints $(\epsilon, 0)$ and $(-\epsilon,0)$ and $\BR$ represents the $x$ direction in $\BHH$ upper half space model. Similarly, $\Sigma^+ = \eta^+ \times\BR$ (See Figure \ref{H-strip}).

Now, consider the hyperbolic isometry $\psi_\lambda(x,y,z)=(\lambda x, \lambda y, \lambda z)$. Let $\psi_\lambda(\Sigma^-)=\Sigma^-_\lambda$, and hence $\Sigma^-_\lambda=\eta^-_\lambda\times \BR$ where $\eta^-_\lambda=\psi_\lambda(\eta^-)$. Let $\lambda_0= \sup\{\lambda \ | \ \eta^-_\lambda \cap\eta^+\neq \emptyset\}$. Clearly, $\psi_1$ is the identity map, and $1\leq \lambda_0<\infty$. However, this implies $\Sigma^-_{\lambda_0}$ and $\Sigma^+$ has tangential intersection as one lies in one side of the other. This contradicts to the maximum principle, Lemma \ref{maxprinciple}.
\end{pf}

\begin{figure}[b]
\begin{center}
$\begin{array}{c@{\hspace{.2in}}c}

\relabelbox  {\epsfxsize=2in \epsfbox{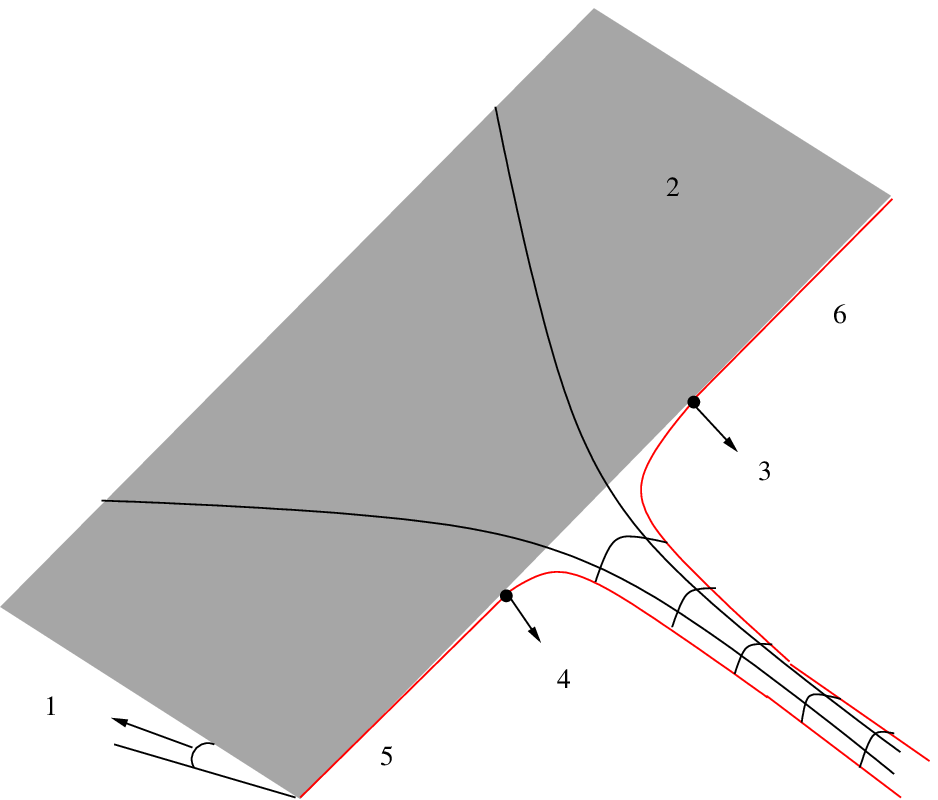}} \relabel{1}{$\theta_H$} \relabel{2}{$\Sigma^\e_H$}  \relabel{3}{$p^-$} \relabel{4}{$p^+$} \relabel{5}{$\beta_\e^+$} \relabel{6}{$\beta_\e^-$} \endrelabelbox &

\relabelbox  {\epsfxsize=1.7in \epsfbox{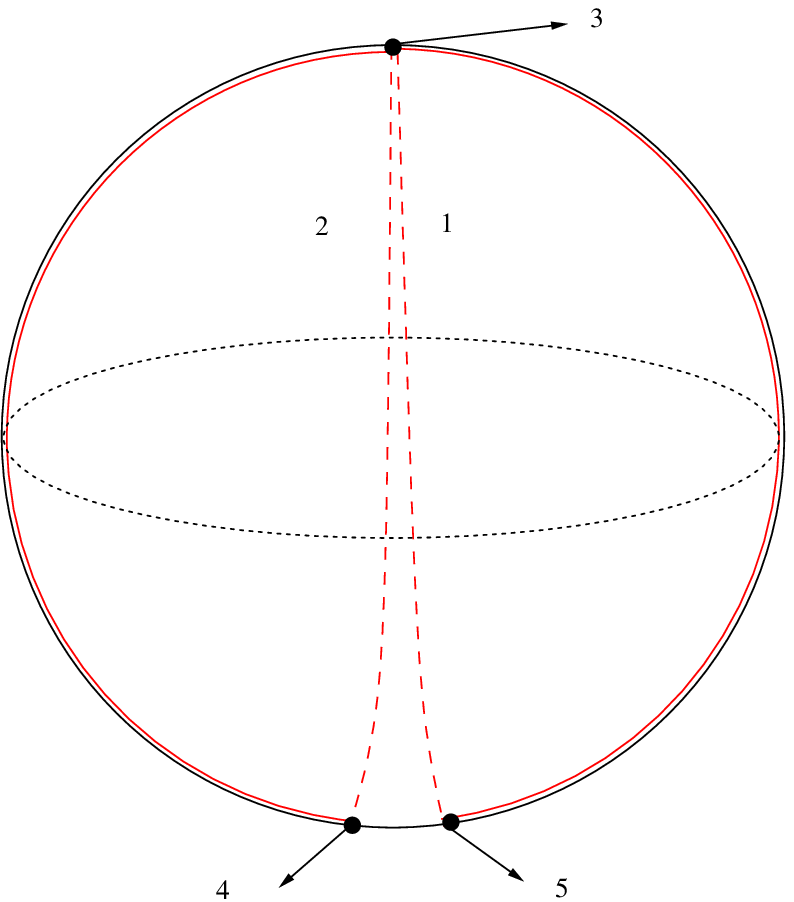}} \relabel{1}{$\beta^+$} \relabel{2}{$\beta^-$} \relabel{3}{$\infty$}  \relabel{4}{$p^-$} \relabel{5}{$p^+$}  \endrelabelbox \\ [0.4cm]

\end{array}$

\end{center}
\caption{\label{H-skillet} \small In the left, $\Sigma^\e_H$ is an $H$-skillet with $\PI \Sigma_H = \beta_\e^+\cup \beta_\e^-$ for $-1<H<1$. In the right, the lines $\beta^+$ and $\beta^-$ are pictured in the Poincare ball model, where $p^\pm$ corresponds to points $(0,\pm 2\e ,0)$. }
\end{figure}

Now, we define $H$-skillets, and show that they are uniquely $H$-minimizing. Again, we use the upper half space model. First, we define its asymptotic boundary $\Gamma$ in $\Si$. Let $u:(-\infty,-1)\cup (1,\infty)\to \BR^+$ be a smooth convex function $u''(x)\geq 0$ such that $u(x)=0$ when $|x|\geq 2$ and $u(x)\to\infty$ when $|x|\to 1$. Define $\Gamma=\mbox{graph}(u)$ in the $xy$-plane, and let $\Omega=\{(x,y) \ | y\leq u(x)\} \cup [-1,1]\times\BR$, i.e. $\partial \Omega = \Gamma$. Similarly, define $\Omega^\e=\psi_\e(\Omega)$, $\Sigma^\e=\psi_\e(\Sigma)$, and $\beta^\pm_\e=\psi_\e(\beta^\pm)$ where $\psi_\e(x,y,z)=(\e x, \e y, \e z)$ is the dilating isometry. Since $\psi_\e$ keep $x$, and $y$ axis fixed, $\Sigma^\e$ is another $H$-skillet with a very thin handle (See Figure \ref{H-skillet}).

It is easy to see that $\Gamma$ is star shaped in $\Si$ with respect to the star point $p^*=(0,-\delta,0)$. Note also that in Poincare ball model, $\Gamma$ looks like union of two star shaped curves $\beta^+$ and $\beta^-$ (with different star points) where they touch each other at one point ($\infty$) (See Figure \ref{H-skillet}).

We claim that $\Gamma$ bounds a unique minimizing $H$-surface $S$, which we call $H$-skillet.\\

\noindent {\bf Claim:} $H$-skillets are uniquely $H$-minimizing.

\vspace{.2cm}

\begin{pf} First notice that $S$ is connected. If $S$ is not connected, then it would bound two symmetric uniquely minimizing $H$-surfaces $\mathcal{P}^+$ and $\mathcal{P}^-$ with $\PI \mathcal{P}^\pm = \beta^\pm$ as they are both star shaped curves \cite{GS}. Similar to $H$-strip case, it is possible to find a spherical $H$-catenoid $\mathcal{C}_H$ with $\PI \mathcal{C}_H= \gamma^+\cup\gamma^-$ where $\gamma^\pm$ are the round circles of radius $c$ with centers $(c+2,c+2,0)$ and $(-c-2,c+2,0)$. Similar to $H$-strip case, $\mathcal{C}_H$ transversely intersect both $\mathcal{P}^+$ and $\mathcal{P}^-$ in simple closed curves. Again for suitable choice of $c$, we get a contradiction as before. The existence of such a $c$ can be seen by using the isometry $\psi_\e$, as $\psi_\e(\Gamma)=\Gamma_\epsilon$ has width $2\e$ instead of $2$ along the skillet handle.

Now, assume that $\Gamma$ bounds more than one minimizing $H$-surface. Then as before, there are canonical minimizing $H$-surfaces $S^-$ and $S^+$ by Lemma \ref{canonical}. Here, we take the exhausting sequence of regions $\{\Omega^-_i\}$ in the side of $p_\e$, i.e. $\bigcup_{i=1}^\infty\Omega^-_i=\Omega$, and for any $i$, $p_\e\subset\Omega_i^-$, and this sequence gives the canonical $S^-$. Consider $S^+$ and $S^-$ in Poincare ball model of $\BHH$. Near the infinity point, the part of $\Gamma$ which corresponds to the skillet edge ($x$-axis) is smooth, say $\tau$. In other words, let $\tau\subset \Gamma \cap N_\e(\infty)$ and $\infty\in \tau$ where $\tau\subset \Gamma$ is the smooth piece. Notice that $\tau$ is an arc in the great circle of $\Si$ corresponding to $x$-axis in upper half space model, i.e. in the upper half space model $\tau=\{|x|>C \mbox{ and } y=0\}\cup \{\infty\}$. Let $T^+=S^+\cap N_\e(\infty)$ and $T^-=S^-\cap N_\e(\infty)$ in the Poincare ball model for some small $\e>0$. By the proof of Lemma \ref{bounreg}, both surfaces $T^+$ and $T^-$ are graphs over $\tau\times [0,\rho)$ for some $\rho>0$. Note that as $S^+$ and $S^-$ are disjoint, so are $T^+$ and $T^-$.

Now, consider the hyperbolic isometry $\psi_\lambda(x,y,z)=(\lambda x, \lambda y, \lambda z)$ again. By construction, for some sufficiently large $\lambda_0$, for any $\lambda\geq \lambda_0$, $\psi_\lambda(S^-)=S^-_\lambda$ would intersect $S^+$. Hence for sufficiently large $\lambda_1$, $T^-_{\lambda_1}\cap T^+$ would be an infinite line $\kappa$ in upper half space model. Hence, in Poincare ball model, $\kappa$ is asymptotic to the point $\infty$, and $\overline{\kappa}$ is a simple closed curve in  $\overline{\BHH}$ with $\infty \in \kappa$. Let $D^-_{\lambda_1} \subset T^-_{\lambda_1}$ and $D^+\subset T^+$ be the $H$-surfaces with boundary $\kappa$. Notice that both $\overline{D}^-_{\lambda_1}$ and $\overline{D}^+$ are both embedded compact disks with boundary $\overline{\kappa}$ in Poincare ball model $\overline{\BHH}$.

Now, we will get a contradiction via maximum principle by moving $D^+$ towards  $D^-_{\lambda_1}$ by isometry. In particular, let $\phi_t$ be the parabolic isometry which fixes the point $0$ (origin in the upper half space model) in $\Si$, and translates $\BHH$ along the great circle in $\Si$ which corresponds to the $y$-axis in upper half space model. Then, $\phi_t(\infty)=q_t$ where $q_t$ is a point in a great circle $\sigma_y$ in $\Si$ corresponding to $y$-axis, i.e. $q_t=(0,-C_t,0)$ in the upper half space model. Here $(0,-C_t,0)$ is the image of $q_t$ in the conversion of Poincare ball model into upper half space model where $t\searrow 0$ implies $C_t\nearrow\infty$. Then, $\phi_t(\tau)=\tau_t$ is an arc in a round circle $\xi_t$ in $\Si$ corresponding to the great circle going through $0$ and $(0,-C_t,0)$, i.e. $\xi_t$ corresponds to $x^2+(y+C/2)^2= C^2/4$ in the upper half space model. This is because $\phi_t(\sigma_x)=\xi_t$ by the definition of the parabolic isometry $\phi_t$.

Let $\phi_t(D^+)=D^+_t$. Then, for sufficiently small $t>0$, $D^+_t\cap D^-_{\lambda_1} \neq \emptyset$ and let $t_1=\sup\{ t \ | \ D^+_t\cap D^-_{\lambda_1} \neq \emptyset\}$. Then, $D^-_{\lambda_1}$ and $D^+_t$ have tangential intersection in an interior point, and one lies in the one side of the other. However, as both $D^-_{\lambda_1}$ and $D^+_t$ are $H$-surfaces, this again contradicts to the maximum principle by Lemma \ref{maxprinciple}.
\end{pf}

\begin{rmk} Notice that $H$-strips and $H$-skillets are defined for $-1<H<1$ instead of $0<H<1$. This is because depending on whether the side we are attaching the bridge is the convex side or concave side of the original surface, the $H$-strips or $H$-skillets can be either positive ($0<H<1$ and mean curvature vector points downwards along the skillet handle) or negative ($-1<H<0$ and mean curvature vector points upwards along the skillet handle). In particular for $0<H<1$, consider the $+H$-skillet $S^+_H$ asymptotic to $\mathcal{P}_H^+ = \{ y=+\frac{H}{\sqrt{1-H^2}} z \}$ and the $-H$-skillet $S^-_H$ asymptotic to  $\mathcal{P}_H^- = \{ y=-\frac{H}{\sqrt{1-H^2}} z\}$. Then, $\PI S^+_H=\PI S^-_H=\Gamma$ define above, and the skillet handles are in the same side ($+y$-axis). However, in $+H$-skillet $S^+_H$, the skillet goes towards the skillet handle, whereas in $-H$-skillet $S^-_H$, the skillet goes away from the skillet handle (See Figure \ref{H-skillet}-left). Similarly, in $+H$-skillet $S^+_H$, the mean curvature vector points downwards along the skillet handle, while in $-H$-skillet $S^-_H$, the mean curvature vector points upwards along the skillet handle.

To see these situations in our constructions in Section 3 and 4, let $\Sigma_1$ be the uniquely minimizing $H$-surface (a spherical cap) with $\PI \Sigma_1 =\Gamma_1$ is a round circle of radius $1$ with center $(0,0,0)$ in upper half space model. Let $\alpha$ be an arc in the unit disk with $\alpha\cap \Gamma_1=\partial \alpha$, and $\alpha\perp \Gamma$. By using Theorem \ref{bridge}, we get a uniquely minimizing $H$-surface $\Sigma_2$. Then, along the bridge the mean curvature vector points upwards, hence the bridge looks like $-H$-strip. In other words, near the endpoints of the bridge $\alpha$, one sees that $\Sigma_2$ looks like $-H$-skillet $S^-_H$. This is true for any bridge $\alpha$ which is in the bounded side of $\Si-\PI\Sigma_1$ in upper half space model.

However, if the endpoints of the bridge $\alpha$ are in different components of $\PI \Sigma$ where $\Sigma$ is a uniquely minimizing $H$-surface, then after applying Theorem \ref{bridge}, we get a uniquely minimizing $H$-surface $\Sigma'$. However this time, along the bridge the mean curvature vector points downwards, hence the bridge looks like $+H$-strip. In other words, near the endpoints of the bridge $\alpha$, one sees that $\Sigma_2$ looks like $+H$-skillet $S^+_H$. Again, this is true for any bridge $\alpha$ which is in the unbounded side of $\Si-\PI\Sigma$ in upper half space model. In particular, the second bridges in the handle cases are examples of this situation.
\end{rmk}

\vspace{-.2cm}

\end{document}